\documentstyle[11pt,leqno,amscd,amssymb,pstricks,latexsym,amsbsy,xypic,mathrsfs,verbatim]{amsart}

\newcommand{\scr}[1]{\mathscr #1}
\newtheorem{thm}{Theorem}[section]
\newtheorem{lem}[thm]{Lemma}
\newtheorem{cor}[thm]{Corollary}
\newtheorem{prop}[thm]{Proposition}

 \theoremstyle{definition}

\newtheorem{rem}[thm]{Remark}
\newtheorem{rems}[thm]{Remarks}

\numberwithin{equation}{thm}

\newcommand{\leb}{\left[}

\newcommand{\rib}{\right]}

\def\lan{\langle}    \def\ran{\rangle}
\def\lr#1{\langle #1\rangle}

\def\az{\alpha}  \def\thz{\theta}
\def\bz{\beta}  \def\dt{\Delta}
\def\gz{\gamma}
  \def\ooz{\Omega}
\def\sz{\sigma}
\def\oz{\omega}
\def\vez{\varepsilon}  
 \def\dz{\delta}

\def\llz{\Lambda}
\def\lz{\lambda}
\def\cM{{\cal M}}      
\def\cC{{\cal C}}    
\def\cH{{\cal H}}     
\def\cZ{{\cal Z}} \def\cO{{\cal O}} 
\def\cI{{\cal I}}

\def\bbN{{\mathbb N}}  \def\bbZ{{\mathbb Z}}  \def\bbQ{{\mathbb Q}}

      \def\vphi{\varphi}
        \def\bfd{{\bf d}} \def\bfU{{\bf U}}
   
  \def\leq{\leqslant}  \def\geq{\geqslant}
 
\def\lra{\longrightarrow}   
\def\lmto{\longmapsto}   
\def\ra{\rightarrow}
\def\Hom{\mbox{\rm Hom}}  
       \def\im{\mbox{\rm Im}\,}
\def\Ext{\mbox{\rm Ext}\,}   
\def\dim{\mbox{\rm dim}\,}   \def\End{\mbox{\rm End}}
\def\Aut{\mbox{\rm Aut}}
\def\udim{{\mathbf dim\,}} \def\oline{\overline}
  
\def\rad{\mbox{\rm rad}\,}  
\def\rep{\mbox{\rm Rep}\,}  \def\soc{\mbox{\rm soc}\,}

\def\fkm{{\frak m}}  
  
  \def\fkp{{\frak p}}

\def\bfb{{\bf b}}
  \def\bfa{{\bf a}}
  \def\bfd{{\bf d}}
\def\tu{{{\widetilde u}}}
\def\dg{{\rm deg}}
\def\blz{{|\lz\rangle}}
\def\bmu{{|\mu\rangle}}
\def\bnu{{|\nu\rangle}}

\def\bempty{{|\emptyset\rangle}}

\def\scrF{{\scr F}}
\def\oln{\overline}
\def\fkM{{\frak M}}
\def\bfcH{\boldsymbol{\cal H}}
\def\bfcC{\boldsymbol{\cal C}}
\def\scrD{\boldsymbol{\cal D}}

\def\wt{\widetilde}
\def\bdc{\boldsymbol{ c}}
\def\bdx{\boldsymbol{ x}}
\def\bdz{\boldsymbol{ z}}

\begin{document}

\title[Hall algebras and $q$-deformed Fock spaces]%
{Hall algebras of cyclic quivers and\\ $q$-deformed Fock spaces}

\author{Bangming Deng and Jie Xiao}
\address{Yau Mathematical Sciences Center, Tsinghua University, Beijing 100084, China}
\email{bmdeng@@math.tsinghua.edu.cn}
\address{Department of Mathematical Sciences,
Tsinghua University, Beijing 100084, China}
\email{jxiao@@math.tsinghua.edu.cn}

\thanks{Supported partially by the Natural Science Foundation of China.}



\begin{abstract} Based on the work of Ringel and Green, one can define the (Drinfeld) double Ringel--Hall algebra
$\scrD(Q)$ of a quiver $Q$ as well as its highest weight modules. The main purpose of the present paper is to show that
the basic representation $L(\llz_0)$ of $\scrD(\dt_n)$ of the cyclic quiver $\dt_n$ provides a realization
of the $q$-deformed Fock space $\bigwedge^\infty$ defined by Hayashi. This is worked out by extending a
construction of Varagnolo and Vasserot. By analysing the structure of nilpotent representations of $\dt_n$,
we obtain a decomposition of the basic representation
$L(\llz_0)$ which induces the Kashiwara--Miwa--Stern decomposition of $\bigwedge^\infty$ and a construction of
the canonical basis of $\bigwedge^\infty$ defined by Leclerc and Thibon in terms of certain monomial basis elements
in $\scrD(\dt_n)$.

\end{abstract}

\maketitle


\section{Introduction}

In \cite{R93}, Ringel introduced the Hall algebra $\bfcH(\dt_n)$ of the cyclic quiver $\dt_n$
with $n$ vertices and showed that its subalgebra generated by simple representations, called
the composition algebra, is isomorphic to the positive part $\bfU^+_v(\widehat{\frak{sl}}_n)$
of the quantized enveloping algebra $\bfU_v(\widehat{\frak{sl}}_n)$. Schiffmann \cite{Sch00} further showed
that $\bfcH(\dt_n)$ is the tensor product of $\bfU_v^+(\widehat{\frak{sl}}_n)$ with
a central subalgebra which is the polynomial ring in infinitely many indeterminates. Following the approach in \cite{X97}, the double Ringel--Hall
algebra $\scrD(\dt_n)$ was defined in \cite{DDF}. Based on \cite{FM,Hub2} and an explicit description
of central elements of $\bfcH(\dt_n)$ in \cite{Hub1}, it was shown in \cite[Th.~2.3.3]{DDF} that $\scrD(\dt_n)$ is
isomorphic to the quantum affine algebra $\bfU_v(\widehat{\frak{gl}}_n)$ defined by Drinfeld's new presentation
\cite{Dr88}.

The $q$-deformed Fock space representation $\bigwedge^\infty$ of the quantized enveloping algebra
$\bfU_v(\widehat{\frak{sl}}_n)$ has been constructed by Hayashi \cite{Hay}, and its crystal basis
was described by Misra and Miwa \cite{MM}. Further, by work of Kashiwara, Miwa, and Stern \cite{KMS},
the action of $\bfU_v(\widehat{\frak{sl}}_n)$ on $\bigwedge^\infty$ is centralized by a Heisenberg algebra
which arises from affine Hecke algebras. This yields a bimodule isomorphism from $\bigwedge^\infty$
to the tensor product of the basic representation of $\bfU_v(\widehat{\frak{sl}}_n)$ and the Fock space
representation of the Heisenberg algebra.

By defining a natural semilinear involution on $\bigwedge^\infty$, Leclerc and Thibon \cite{LT}
obtained in an elementary way a canonical basis of $\bigwedge^\infty$. It was conjectured in \cite{LLT,LT}
that for $q=1$, the coefficients of the transition matrix of the canonical basis on the natural
basis of $\bigwedge^\infty$ are equal to the decomposition numbers for Hecke algebras and quantum Schur algebras
at roots of unity. These conjecture have been proved, respectively, by Ariki \cite{A0} and
Varagnolo and Vasserot \cite{VV}. For the categorification of the Fock space, see, for example,
\cite{Shan,HY,SW}.

In \cite{VV}, Varagnolo and Vasserot extended the $\bfU_v(\widehat{\frak{sl}}_n)$-action on
the Fock space $\bigwedge^\infty$ to that of the extended Ringel--Hall algebra $\scrD(\dt_n)^{\leq 0}$
of the cyclic quiver $\dt_n$. They also showed that the canonical basis of the Ringel--Hall algebra $\bfcH(\dt_n)$
in the sense of Lusztig induces a basis of $\bigwedge^\infty$ which conjecturally coincides with the
canonical basis constructed by Leclerc--Thibon \cite{LT}. This conjecture was proved by
Schiffmann \cite{Sch00} by identifying the central subalgebra of $\bfcH(\dt_n)$ with the ring of symmetric functions.

The main purpose of the present paper is to extend Varagnolo--Vasserot's construction to obtain a
$\scrD(\dt_n)$-module structure on the Fock space $\bigwedge^\infty$ which is shown to be
isomorphic to the basic representation $L(\llz_0)$ of $\scrD(\dt_n)$. Moreover, the central
elements in the positive and negative parts of $\scrD(\dt_n)$ constructed by Hubery \cite{Hub1} give rise naturally to
the operators introduced in \cite{KMS} which generate the Heisenberg algebra. Furthermore, the structure of
$\scrD(\dt_n)$ yields a decomposition of $L(\llz_0)$ which induces
the Kashiwara--Miwa--Stern decomposition of $\bigwedge^\infty$. This also provides a way to construct
the canonical basis of $\bigwedge^\infty$ in \cite{LT} in terms of
certain monomial basis elements of $\scrD(\dt_n)$.

The paper is organized as follows. In Section 2 we review the classification of (nilpotent)
representations of both infinite linear quiver $\dt_\infty$ and the cyclic quiver $\dt_n$
with $n$ vertices and discuss their generic extensions. Section 3 recalls
the definition of Ringel--Hall algebras $\bfcH(\dt_\infty)$ and $\bfcH(\dt_n)$ of $\dt_\infty$ and
$\dt_n$ as well as the maps from the homogeneous spaces of $\bfcH(\dt_n)$ to those of $\bfcH(\dt_\infty)$
introduced in \cite{VV}. The images of basis elements of $\bfcH(\dt_n)$ under these maps are described.
In Section 4 we first
follow the approach in \cite{X97} to present the construction of double Ringel--Hall algebras
of both $\dt_\infty$ and $\dt_n$ and then study the irreducible highest weight $\scrD(\dt_n)$-modules
based on the results in \cite{JKK}. Section 5 recalls from \cite{Hay,MM,VV} the Fock space representation
$\bigwedge^\infty$ over $\bfU_v(\widehat{\frak{sl}}_\infty)$ ($\cong\scrD(\dt_\infty)$) as well as over
$\bfU^+_v(\widehat{\frak{sl}}_n)$. In Section 6 we define the $\scrD(\dt_n)$-module structure on $\bigwedge^\infty$ based on
\cite{KMS,VV}. It is shown in Section 7 that $\bigwedge^\infty$
is isomorphic to the basic representation of $\scrD(\dt_n)$. In the final section,
we present a way to construct the canonical basis of $\bigwedge^\infty$ and interpret the ``ladder method''
construction of certain basis elements in $\bigwedge^\infty$ in terms of generic extensions of
nilpotent representations of $\dt_n$.

\section{Nilpotent representations and generic extensions}

In this section we consider nilpotent representations of both a cyclic quiver $\dt=\dt_n$ with $n$ vertices
($n\geq 2$) and the infinite quiver $\dt=\dt_\infty$ of type $A_\infty^\infty$ and study their generic extensions.
We show that the degeneration order of nilpotent representations of $\dt_n$ induces the dominant order
of partitions.

Let $\dt_\infty$ denote the infinite quiver of type $A^\infty_\infty$
\begin{center}
\begin{pspicture}(-3,-.4)(3.6,0.4)
\psset{linewidth=0.5pt, arrowsize=3.5pt}
\psset{xunit=.8cm,yunit=.7cm}
\psdot*(-3.9,0) \psdot*(-2.6,0) \psdot*(-1.3,0) \psdot*(0,0)
\psdot*(1.3,0)\psdot*(2.6,0) \psdot*(3.9,0)
\uput[d](-2.6,0){$_{-2}$} \uput[d](-1.3,0){$_{-1}$} \uput[d](0,0){$_0$}
\uput[d](1.3,0){$_1$} \uput[d](2.6,0){$_2$}
\psline{->}(-3.9,0)(-2.6,0) \psline{->}(-2.6,0)(-1.3,0)
\psline(-1.3,0)(0,0) \psline{->}(0,0)(1.3,0) \psline{->}(1.3,0)(2.6,0)
\psline{->}(2.6,0)(3.9,0)
\psline[linestyle=dotted,linewidth=1pt](-5.2,0)(-3.9,0)
\psline[linestyle=dotted,linewidth=1pt](3.9,0)(5.2,0)
\end{pspicture}
\end{center}
 with vertex set $I=I_\infty=\bbZ$, and for $n\geq 2$, let $\dt_n$ denote the cyclic quiver
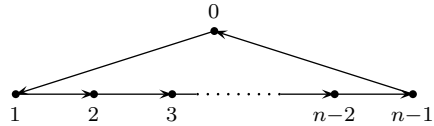
\begin{center}
\begin{pspicture}(-3,-.6)(3.6,1.6)
\psset{linewidth=0.5pt, arrowsize=3.5pt}
\psset{xunit=.8cm,yunit=.7cm}
\psdot*(-3,0) \psdot*(-1.7,0) \psdot*(-.4,0) \psdot*(2.3,0)
\psdot*(3.6,0)\psdot*(.3,1.2) \uput[u](.3,1.2){$_{0}$}
\uput[d](-3,0){$_1$} \uput[d](-1.7,0){$_2$} \uput[d](-.4,0){$_3$}
\uput[d](2.3,0){$_{n-2}$} \uput[d](3.6,0){$_{n-1}$}
\psline{->}(-3,0)(-1.7,0) \psline{->}(-1.7,0)(-.4,0)
\psline(-.4,0)(0,0)
\psline[linestyle=dotted,linewidth=1pt](0,0)(1.4,0)\psline{->}(1.4,0)(2.3,0)
\psline{->}(2.3,0)(3.6,0) \psline{->}(3.6,0)(.3,1.2)
\psline{->}(.3,1.2)(-3,0)
\end{pspicture}
\end{center}

\noindent with vertex set $I=I_n=\bbZ/n\bbZ=\{0,1,\ldots,n-1\}$. For each $i\in I_\infty=\bbZ$, let $\bar i$
denote its residue class in $I_n=\bbZ/n\bbZ$. We also simply write $\bar i\pm 1$ to denote the
residue class of $i\pm1$ in $\bbZ/n\bbZ$.

Given a field $k$, we denote by $\rep^0\dt$ the category of finite dimensional nilpotent representations
of $\dt$ ($=\dt_\infty$ or $\dt_n$) over $k$. (Note that each finite dimensional representation of $\dt_\infty$ is
automatically nilpotent.) Given a representation $V=(V_i,V_\rho)\in\rep^0\dt$, the vector $\udim V=(\dim_k V_i)_{i\in I}$
is called the {\it dimension vector} of $V$. The Grothendieck group of $\rep^0\dt$ is identified with
the free abelian group $\bbZ I$ with basis $I$. Let $\{\vez_i\mid i\in I\}$ denote the standard basis
of $\bbZ I$. Thus, elements in $\bbZ I$ will be written as $\bfd=(d_i)_{i\in I}$ or
$\bfd=\sum_{i\in I} d_i \vez_i$. In case $I=\bbZ/n\bbZ$, we sometimes write $\bbZ^n$ for $\bbZ I$.

 The Euler form $\langle-,-\rangle:\bbZ I\times \bbZ I\ra\bbZ$ is defined by
 $$\lr{\udim M, \udim N}=\dim_k\Hom_{k\dt}(M, N)-\dim_k\Ext^1_{k\dt}(M, N).$$
Its symmetrization
$$(\udim M, \udim N)=\lr{\udim M, \udim N} +\lr{\udim N, \udim M}$$
 is called the symmetric Euler form.

It is well known that the isoclasses of representations in $\rep^0\dt$ are parametrized
by the set $\frak M$ consisting of all multisegments
$$\fkm=\sum_{i\in I,\,l\geq 1} m_{i,l}[i,l),$$
 where all $m_{i,l}\in\bbN$, but finitely many, are zero. More precisely, the representation
$M(\fkm)=M_k(\fkm)$ associated with $\fkm$ is defined by
$$M(\fkm)=\bigoplus_{i\in I,l\geq 1}m_{i,l} S_i[l],$$
 where $S_i[l]$ denotes the representation of $\dt$ with the simple top $S_i$ and length $l$.
For each $\bfd\in\bbN I$, put
$$\fkM^\bfd=\{\fkm\in\fkM\mid \udim M(\fkm)=\bfd\}.$$
 Furthermore, we will write $\fkM=\fkM_\infty$ (resp., $\fkM=\fkM_n$) if $I=\bbZ$ (resp., $I=\bbZ/n\bbZ$).

It is also known that there exist Auslander--Reiten sequences in $\rep^0\dt$, that is, for each $M\in\rep^0\dt$,
there is an Auslander--Reiten sequence
$$0\lra\tau M\lra E\lra M\lra 0,$$
where $\tau M$ denotes the Auslander--Reiten
translation of $M$. It is clear that $\tau$ induces an isomorphism $\tau:\bbZ I\ra \bbZ I$
such that $\tau(\udim M)=\udim\tau M$. In particular, $\tau(\vez_i)=\vez_{i+1},\;\forall\,i\in I$.
If $\dt=\dt_n$, then $\tau^{sn}={\rm id}$ for all $s\in\bbZ$. For $\fkm\in\fkM$, let $\tau\fkm$ be
defined by $M(\tau\fkm)\cong \tau M(\fkm)$.

Given $\bfd\in\bbN I$, let $V=\oplus_{i\in I}V_i$ be an $I$-graded vector space with
dimension vector $\bfd$. Consider
$$E_V=\{(x_i)\in\bigoplus_{i\in I}\Hom_k(V_i,V_{i+1})\mid \text{$x_{n-1}\cdots x_0$ is nilpotent if $\dt=\dt_n$.}\}.$$
 Then each element $x\in E_V$ defines a representation $(V,x)$ of dimension vector $\bfd$ in $\rep^0\dt$.
Moreover, the group
$$G_V=\prod_{i\in I}{\rm GL}(V_i)$$
 acts on $E_V$ by conjugation, and there is a bijection between the $G_V$-orbits and the isoclasses
of representations in $\rep^0\dt$ of dimension vector $\bfd$. For each $x\in E_V$, by $\cO_x$ we denote
the $G_V$-orbit of $x$. In case $k$ is algebraically closed, we have the equalities
\begin{equation}\label{dim-orbit}
\dim\cO_x=\dim G_V-\dim \End_{k\dt}(V,x)=\sum_{i\in I}d_i^2-\dim \End_{k\dt}(V,x).
\end{equation}
By abuse of notation, for each $M\in\rep^0\dt$, we denote by $\cO_M$ the orbit of $M$.

Following \cite{Bon, Rein, DD1}, given two representations $M,N$ in $\rep^0\dt$,
there exists a unique (up to isomorphism) extension $G$ of $M$ by
$N$ such that $\dim\End_{k\dt}(G)$ is minimal.
The extension $G$ is called the {\it generic extension} of $M$ by
$N$, denoted by $M*N$. Moreover, generic extensions satisfy the associativity, i.e., for $L,M,N\in\rep^0\dt$,
$$L\ast (M\ast N)\cong (L\ast M)\ast N.$$
 Let $\cM(\dt)$ denote the the set of isoclasses of representations in $\rep^0\dt$. Define a multiplication
on $\cM(\dt)$ by setting
$$[M]\ast [N]=[M\ast N].$$
 Then $\cM(\dt)$ is a monoid with identity $[0]$, the isoclass of zero representation of $\dt$.

 By \cite{Rein,DD1}, the generic extension $M\ast N$ can be also characterized as the unique maximal element
among all the extensions of $M$ by $N$ with respect to the degeneration
order $\leq_\dg$ which is defined by setting $M\leq_\dg N$ if $\udim M=\udim N$ and
\begin{equation}\label{deg-order-Hom}
\dim_k\Hom_{k\dt}(M,X)\geq\dim_k\Hom_{k\dt}(N,X),\;\,\text{ for all $X\in\rep^0\dt$.}
\end{equation}
 If $k$ is algebraically closed, then $M\leq_\dg N$ if and only if $\overline{\cO}_M\subseteq \cO_N$,
where $\overline{\cO}_M$ is the closure of $\cO_M$.
 This defines a partial order relation on the set $\cM(\dt)$ of isoclasses of representations in $\rep^0\dt$;
see \cite[Th.~2]{Zwa} or \cite[Lem.~3.2]{DD1}. By \cite[2.4]{Rein}, for $M,N,M',N'\in\rep^0\dt$,
$$M'\leq_\dg M,N'\leq_\dg N\Longrightarrow M'*N'\leq_\dg M*N.$$
For $\fkm,\fkm'\in\fkM_n$ (resp., $\fkM_\infty$), we write $\fkm\leq_\dg \fkm'$ (resp., $\fkm\leq^\infty_\dg \fkm'$) if
$M(\fkm)\leq_\dg M(\fkm')$ in $\rep^0\dt_n$ (resp., $\rep\dt_\infty$).

By \cite{BG82,Gab81}, there is a covering functor
$${\scr F}:\rep \dt_\infty\lra \rep^0 \dt_n$$
 sending $S_i[l]$ to $S_{\bar i}[l]$ for $i\in\bbZ$ and $l\geq 1$. Moreover, $\scrF$ is dense and exact,
and the Galois group of $\scrF$ is the infinite cyclic group $G$ generated by $\tau^n$, i.e., $\tau^n(S_i[l]=S_{i+n}[l])$.
 For $\fkm\in\fkM_\infty$, let $\scrF(\fkm)\in\fkM_n$ be such that
$M(\scrF(\fkm))\cong \scrF(M(\fkm))\in\rep^0 \dt_n$. From \eqref{deg-order-Hom} we easily deduce that
for $M,N\in\rep \dt_\infty$,
\begin{equation}\label{covering-preserving-dg}
M\leq_\dg N\Longrightarrow\scrF(M)\leq_\dg \scrF(N).
\end{equation}

%

%
%
%

The following two classes of representations will play an important role later on.
For each $\bfd=(d_i)\in\bbN I$, we set
$$ S_\bfd=\bigoplus_{i\in I}d_i S_i[1]\in\rep^0\dt.$$
 In other words, $S_\bfd$ is the unique semisimple representation of dimension vector $\bfd$.

Let $\Pi$ be the set of all partitions $\lz=(\lz_1,\ldots,\lz_t)$ (i.e.,
$\lz_1\geq\cdots\geq \lz_t\geq 1$). For each $\lz\in\Pi$, define
$$\fkm_\lz=\sum_{s=1}^t[1-s,\lz_s)\in\fkM.$$
  Then
$$M(\fkm_\lz)=S_0[\lz_1]\oplus S_{-1}[\lz_2]\oplus\cdots\oplus S_{1-t}[\lz_t]\in\rep^0\dt.$$
 If $\dt=\dt_\infty$, then we sometimes write $\fkm_\lz=\fkm^\infty_\lz\in\fkM_\infty$ to make a distinction.
 It follows from the definition that
$\scrF(\fkm_\lz^\infty)=\fkm_\lz$ for all $\lz\in\Pi$.

\begin{prop} \label{order-coincidence} Let $\lz,\mu\in\Pi$.

\begin{itemize}
\item[(1)] If $\dt=\dt_\infty$, then
$$\udim M(\fkm_\mu^\infty)=\udim M(\fkm_\lz^\infty)\Longleftrightarrow \mu=\lz.$$
 In particular, for each $\fkm\in\fkM_\infty$, there exists at most one $\nu\in\Pi$
such that $\fkm=\fkm^\infty_\nu$.

\item[(2)] If $\dt=\dt_n$, then
$$M(\fkm_\mu)\leq_\dg M(\fkm_\lz)\Longrightarrow \mu\trianglelefteq \lz,$$
 where $\trianglelefteq$ is the dominance order on $\Pi$, i.e., $\mu \trianglelefteq \lz\Longleftrightarrow
 \sum_{j=1}^i\mu_j\leq \sum_{j=1}^i\lz_j$, $\forall\,i\geq 1$.
\end{itemize}
\end{prop}

\begin{pf} (1) By definition, both the socles of $M(\fkm_\lz^\infty)$ and $M(\fkm_\mu^\infty)$ are
multiplicity-free. Thus, comparing the socles of $S_0[\lz_1]$ and $S_0[\mu_1]$ gives $\lz_1=\mu_1$.
The lemma then follows from an inductive argument.

(2) Suppose $M(\fkm_\mu)\leq_\dg M(\fkm_\lz)$. By viewing $\fkm_\lz$ and
$\fkm_\mu$ as multipartitions in $\fkM_n$, we obtain by \cite[Prop.~2.7]{DDM} that for each $l\geq 1$,
$$\sum_{s=1}^l\wt\mu_s\geq \sum_{s=1}^l \wt\lz_s,$$
 where $\wt\lz=(\wt \lz_1,\wt\lz_2,\ldots)$ and $\wt\mu=(\wt \mu_1,\wt\mu_2,\ldots)$ are the dual
partition of $\lz$ and $\mu$, respectively, that is, $\wt\mu\unrhd \wt\lz$. By \cite[1.1]{Mac},
 $\mu\trianglelefteq \lz$.
\end{pf}

\section{Ringel--Hall algebra of the quiver $\dt$}

In this section we introduce the Ringel--Hall algebra $\cH(\dt)$ of $\dt$ ($=\dt_n$ or $\dt_\infty$) and the maps
from homogeneous subspaces of $\cH(\dt_n)$ to those of $\cH(\dt_\infty)$ defined in \cite[6.1]{VV}. We also
describe the images of basis elements of $\cH(\dt_n)$ under these maps.

The cyclic quiver $\dt_n$ gives the $n\times n$ Cartan matrix $C_n=(a_{ij})_{i,j\in I}$ of type
$\widehat A_{n-1}$,
while $\dt_\infty$ defines the infinite Cartan matrix $C_\infty=(a_{ij})_{i,j\in\bbZ}$.
 Thus, we have the associated quantum enveloping algebras $\bfU_v(\widehat{\frak{sl}}_n)$ and
$\bfU_v({\frak{sl}}_\infty)$ which are $\bbQ(v)$-algebras with generators
$K_i^{\pm1},E_i,F_i, D$ ($i\in I=\bbZ/\bbZ_n$) and $K_i^{\pm 1}, E_i,F_i$ ($i\in \bbZ$),
respectively, and the quantum Serre relations. In particular, the relations involving the
generator $D$ in $\bfU_v(\widehat{\frak{sl}}_n)$ are
$$D D^{-1}=1=D^{-1}D,\,K_i D=DK_i,\,DE_i=v^{\dz_{0,i}}E_i D,\,DF_i=v^{-\dz_{0,i}}F_i D,\;\;\forall\,i\in I;$$
 see \cite[Def.~3.16]{Ariki}. The subalgebra of $\bfU_v(\widehat{\frak{sl}}_n)$ generated
by $K_i^{\pm1},E_i,F_i$ ($i\in I=\bbZ/\bbZ_n$) is denoted by $\bfU'_v(\widehat{\frak{sl}}_n)$.

By \cite{R90a, R93, Guo}, for $\fkp,\fkm_1,\ldots,\fkm_t\in\fkM$,
there is a polynomial $\vphi^{\fkp}_{\fkm_1,\ldots, \fkm_t}(q)\in\bbZ[q]$ (called Hall polynomial)
such that for each finite field $k$,
$$\vphi^{\fkp}_{\fkm_1,\ldots,\fkm_t}(|k|)=F^{M_k(\fkp)}_{M_k(\fkm_1),\ldots, M_k(\fkm_t)},$$
 which is by definition the number of the filtrations
$$M=M_0\supseteq M_1\supseteq \cdots \supseteq M_{t-1}\supseteq M_t=0$$
such that $M_{s-1}/M_s\cong M_k(\fkm_s)$ for all $1\le s \le t$.
 It is also known that for each $\fkm\in\fkM$, there is a polynomial $a_\fkm(q)\in\bbZ[q]$ such that
for each finite field $k$,
$$a_\fkm(|k|)=|\Aut_{k\dt}(M_k(\fkm))|.$$

Let $\cZ=\bbZ[v,v^{-1}]$ be the Laurent polynomial ring over $\bbZ$ in indeterminate $v$. By definition,
the (twisted generic) {\it Ringel--Hall algebra} $\cH(\dt)$ of $\dt$ is the free $\cZ$-module with basis
$\{u_{\fkm}|\fkm\in\fkM\}$ and multiplication given by
\begin{equation}\label{Hall-multiplication}
u_\fkm  u_{\fkm'}=v^{\lr{\udim M(\fkm),\udim M(\fkm')}}\sum_{\fkp\in\fkM}\vphi^\fkp_{\fkm,\fkm'}(v^2)u_\fkp.
\end{equation}
In practice, we also write $u_{\fkm}=u_{[M(\fkm)]}$ in order to
make certain calculations in terms of modules. Furthermore, for each $\bfd\in\bbN I$, we simply
write $u_\bfd=u_{[S_\bfd]}$. Moreover, both $\cH(\dt)$ and $\cC(\dt)$ are $\bbN I$-graded:
\begin{equation}\label{grading}
\cH(\dt)=\oplus_{\bfd\in\bbN I}\cH(\dt)_\bfd\;\text{ and }\;\cC(\dt)=\oplus_{\bfd\in\bbN I}\cC(\dt)_\bfd,
\end{equation}
where $\cH(\dt)_\bfd$ is spanned by all $u_\fkm$ with $\fkm\in\fkM^\bfd$ and
$\cC(\dt)_\bfd=\cC(\dt)\cap\cH(\dt)_\bfd$. Since the Auslander--Reiten translate $\tau:\rep^0\dt\ra \rep^0\dt$ is an auto-equivalence,
it induces an automorphism
$\tau:\cH(\dt)\ra \cH(\dt),u_\fkm\mapsto u_{\tau\fkm}$. We also consider the $\bbQ(v)$-algebra
$$\bfcH(\dt)=\cH(\dt)\otimes_\cZ\bbQ(v).$$

\begin{rem} We remark that the Hall algebra of $\dt$ defined in \cite{VV} is the opposite algebra
of $\cH(\dt)$ given here with $v$ being replaced by $v^{-1}$. Thus, $v$ and $v^{-1}$ should be
swaped when comparing with the formulas in \cite{VV}.
\end{rem}


For each $i\in I$, set $u_i=u_{[S_i]}$. We then denote by
$\cC(\dt)$ the subalgebra of $\cH(\dt)$ generated by the divided power
$u_i^{(t)}=u_i^t/[t]^!$, $i\in I$ and $t\geq 1$, called the {\it composition algebra} of
$\dt$, where $[t]^!=[t][t-1]\cdots [1]$ with $[m]=(v^m-v^{-m})/(v-v^{-1})$. It is known that
$\cC(\dt_\infty)=\cH(\dt_\infty)$ and there is an isomorphism
$\bfU^+_v({\frak{sl}}_\infty)\cong \bfcH(\dt_\infty)$ taking $E_i\mapsto u_i,\;\forall\,i\in I_\infty=\bbZ$.
 But, for $n\geq 2$, $\cC(\dt_n)$ is a proper subalgebra of $\cH(\dt_n)$. By \cite{R93},
$$\bfU^+_v(\widehat{\frak{sl}}_n)\cong\bfcC(\dt_n):=\cC(\dt_n)\otimes_\cZ\bbQ(v),\;E_i\longmapsto u_i,\;\forall\,i\in I_n.$$
 By \cite[Th.~2.2]{Sch00}, $\bfcH(\dt_n)$ is decomposed into the tensor product of $\bfcC(\dt_n)$
and a polynomial ring in infinitely many indeterminates which are central elements in $\cH(\dt_n)$.
Such central elements have been explicitly constructed in \cite{Hub1}.
More precisely, for each $t\geq 1$, let
\begin{equation}\label{central elements}
\bdc_t=(-1)^tv^{-2nt}\sum_{\fkm}(-1)^{{\rm dim}\,{\rm End}(M(\fkm))}a_\fkm(v^2) u_\fkm \in\cH(\dt_n),
\end{equation}
 where the sum is taken over all $\fkm\in \fkM_n$ such that $\udim M(\fkm)=t\dz$
 with $\dz=(1,\ldots,1)\in\bbN I_n$, and
${\rm soc}\,M(\fkm)$ is square-free, i.e., $\udim\, {\rm
soc}\,M(\fkm)\leq \dz$. The following result is proved in \cite{Hub1}.

\begin{thm} \label{Schiffmann-Hubery}  The elements $c_m$ are central in $\cH(\dt_n)$.
Moreover, there is a decomposition
$$\bfcH(\dt_n)= \bfcC(\dt_n) \otimes_{\bbQ(v)}{\bbQ}(v)[\bdc_1,\bdc_2,\ldots],$$
 where ${\bbQ}(v)[\bdc_1,\bdc_2,\ldots]$ is the polynomial algebra in
$\bdc_t$ for $t\geq 1$. In particular, $\bfcH(\dt_n)$ is generated by
$u_i$ and $\bdc_t$ for $i\in I_n$ and $t\geq 1$.
\end{thm}

For each $\fkm\in\fkM$, set
$d(\fkm)=\dim M(\fkm)$, $\bfd(\fkm)=\udim M(\fkm)$ and define
\begin{equation} \label{tilde-u}
\wt u_\fkm=v^{{\rm dim\,End}_{k\dt}(M(\fkm))-d(\fkm)}u_\fkm.
\end{equation}
 Then $\{\wt u_\fkm\mid \fkm\in\fkM\}$ is also a $\cZ$-basis of $\cH(\dt)$ which plays a
role in the construction of the canonical basis. In particular,
$$\wt u_i=u_i\;\text{ for each $i\in I$ and }\; \wt u_\bfd=v^{\sum_i(d^2_i-d_i)}u_\bfd\;
\text{ for each $\bfd\in\bbN I$.}$$

Consider the map $\pi:\bbZ I_\infty\ra \bbZ I_n,\bfd\mapsto\bar\bfd$, where
$\pi(\bfd)=\bar\bfd=(d_{\bar i})$ is defined by
$$d_{\bar i}=\sum_{j\in \bar i}d_j,\;\, \forall\,\bar i\in I_n=\bbZ/n\bbZ.$$
 Then for each representation $M\in \rep \dt_\infty$,
$$\udim \scrF(M)=\pi(\udim M).$$

Take $\bfd\in \bbN I_\infty$ with $\bar\bfd=\pi(\bfd)$. By \cite[6.1]{VV}, there is
a $\cZ$-linear map
$$\gz_\bfd: \cH(\dt_n)_{\bar\bfd}\lra \cH(\dt_\infty)_{\bfd}.$$
  The first two statements in the following lemma are taken from \cite[Sect.~6.1]{VV}, and the third one
follows from the isomorphism $\tau:\cH(\dt_\infty)\ra \cH(\dt_\infty)$.

\begin{lem} \label{compatibility-gamma} {\rm (1)} For each $\bfd\in \bbN I_\infty$,
$\gz_\bfd(\wt u_{\bar\bfd})=v^{-h(\bfd)}\wt u_\bfd$, where $h(\bfd)=\sum_{i<j,\bar i=\bar j}d_i(d_{j+1}-d_j)$.

{\rm (2)} Fix $\az,\beta\in \bbN I_n$ with $\bar\bfd=\az+\beta$. Then for $x\in\cH(\dt_n)_\az$ and
$y\in\cH(\dt_n)_\beta$,
\begin{equation}\label{compat-multi-gamma}
\sum_{\bfa,\bfb}v^{\kappa(\bfa,\bfb)}\gz_\bfa(x)\gz_\bfb(y)=\gz_\bfd(xy),
\end{equation}
 where the sum is taken over all pairs $\bfa,\bfb\in\bbN I_\infty$ satisfying
$\bfa+\bfb=\bfd$, $\bar\bfa=\az$, and $\bar\bfb=\beta$, and
$\kappa(\bfa,\bfb)=\sum_{i>j,\bar i=\bar j}a_i(2b_j-b_{j-1}-b_{j+1}).$

{\rm (3)} For each $\bfd\in \bbN I_\infty$ and $\fkm\in\fkM_n^{\bar \bfd}$, $\gz_{\tau^n(\bfd)}(\wt u_\fkm)
=\tau^n(\gz_\bfd(\wt u_\fkm))$.
\end{lem}

We now describe the images of the basis elements of $\cH(\dt_n)_{\bar\bfd}$ under $\gz_\bfd$.

\begin{prop} \label{images-under-gamma} Let $\bfd\in\bbN I_\infty$ and $\fkm\in\fkM_n$ be
such that $\az:=\udim M(\fkm)=\oln\bfd$. Then
$$\gz_\bfd(\wt u_\fkm)\in\sum_{{\frak z}\in\fkM_\infty,\,\scrF(\frak z)\leq_\dg \fkm}\cZ \wt u_{\frak z}.$$
\end{prop}

\begin{pf} Consider the radical filtration of $M=M(\fkm)$
$$M=\rad^0 M\supseteq \rad M \supseteq\cdots\supseteq \rad^{\ell-1}M\supseteq \rad^\ell M=0$$
 with $\rad^{s-1}M/\rad^s M\cong S_{\az_s}$, where $\ell$ is the Loewy length of of $M$ and $\az_s\in\bbN I_n$ for $1\leq s\leq \ell$.
 Then $M=S_{\az_1}\ast\cdots \ast S_{\az_\ell}$. Moreover, by \cite[Sect.~9]{DDX},
$$\wt u_{\az_1}\cdots \wt u_{\az_\ell}=\wt u_\fkm+\sum_{\fkp<_\dg \fkm}f_{\fkm,\fkp}\wt u_\fkp,\;
\text{ where $f_{\fkm,\fkp}\in\cZ$.}$$

On the one hand, by induction with respect to the order $\leq_\dg$, we may assume that
for each $\fkp\in\fkM_n^\bfd$ with $\fkp<_\dg \fkm$, $\gz_\bfd(\wt u_\fkp)$ is a $\cZ$-linear combination of $\wt u_{\frak y}$
with ${\frak y}\in\fkM_\infty$ satisfying $\scrF(\frak y)\leq_\dg \fkp$. Therefore,
$$\gz_\bfd(\wt u_\fkm)=\gz_\bfd(\wt u_{\az_1}\cdots \wt u_{\az_\ell})+x,$$
  where $x=-\sum_{\fkp<_\dg \fkm}f_{\fkm,\fkp}\gz_\bfd(\wt u_\fkp)$ is a $\cZ$-linear combination of
$\wt u_{\frak z}$ with $\scrF(\frak z)<_\dg \fkm$.

On the other hand, by applying \eqref{compat-multi-gamma} inductively, we obtain
$$\gz_\bfd(\wt u_{\az_1}\cdots \wt u_{\az_\ell})=
\sum_{\bfa_1,\ldots,\bfa_\ell}v^{\sum_{s<t}\kappa(\bfa_s, \bfa_t)-\sum_s h(\bfa_s)}
\wt u_{\bfa_1}\cdots \wt u_{\bfa_\ell},$$
 where the sum is taken over all sequences $\bfa_1,\ldots,\bfa_\ell\in\bbN I_\infty$
satisfying
$$\bfa_1+\cdots+\bfa_\ell=\bfd\;\text{ and }\; \oln{\bfa_s}=\az_s,\;\forall\, 1\leq s\leq\ell.$$
  By the definition, each term $\wt u_{\bfa_1}\cdots \wt u_{\bfa_\ell}$ is a $\cZ$-linear
combination of $\wt u_{\frak y}$ such that $M(\frak y)$ admits a filtration
$$M({\frak y})=X_0\supset X_1\supset \cdots \supset X_{\ell-1}\supset X_\ell=0$$
 such that $X_{s-1}/X_s\cong S_{\bfa_s}$ for all $1\leq s\leq\ell$. Applying the exact functor $\scrF$ gives a
 filtration of $\scrF(M(\frak y))$
$$\scrF(M({\frak y}))=\scrF(X_0\supset \scrF(X_1))\supset \cdots \supset \scrF(X_{\ell-1})\supset \scrF(X_\ell)=0$$
 such that
$$\scrF(X_{s-1})/\scrF(X_s)\cong \scrF(X_{s-1}/X_s)\cong S_{\az_s},\;\forall\,1\leq s\leq\ell.$$
 Therefore,
$$\scrF(M({\frak y}))=M(\scrF(\pi))\leq_\dg S_{\az_1}\ast\cdots \ast S_{\az_\ell}=M(\fkm),$$
 that is, $\scrF({\frak y})\leq_\dg \fkm$.

In conclusion, we obtain that
$$\gz_\bfd(\wt u_\fkm)\in\sum_{{\frak z}\in\fkM_\infty,\,\scrF(\frak z)\leq_\dg \fkm}\cZ \wt u_{\frak z}.$$
\end{pf}

Fix $\lz\in\Pi$ and write
$$\bfd(\lz)=\udim M(\fkm_\lz^\infty)\in\bbN I_\infty\;\text{ and }\; \az(\lz)=\udim M(\fkm_\lz)\in\bbN I_n.$$
 By the definition of $M(\fkm_\lz^\infty)$ and $M(\fkm_\lz)$, the radical filtration of $\wt M=M(\fkm_\lz^\infty)$
$$\wt M=\rad^0 \wt M\supseteq \rad \wt M \supseteq\cdots\supseteq \rad^{\ell-1}\wt M\supseteq \rad^\ell \wt M=0$$
 gives rise to the radical filtration of $M(\fkm_\lz)=\scrF(\wt M)$
$$M(\fkm_\lz)=\scrF(\rad^0 \wt M)\supseteq \scrF(\rad \wt M)\supseteq\cdots\supseteq \scrF(\rad^{\ell-1}\wt M)\supseteq \scrF(\rad^\ell \wt M)=0,$$
 that is, $\scrF(\rad^s \wt M)=\rad^s(M(\fkm_\lz))$ for $1\leq s\leq \ell$.
  Let $\bfd(\lz)_s\in\bbN I_\infty$ and $\az(\lz)_s\in\bbN I_n$, $1\leq s\leq \ell$, be such that
$$\rad^{s-1} \wt M=\rad^s \wt M\cong S_{\bfd(\lz)_s}\text{ and }\;\rad^{s-1}M(\fkm_\lz)/\rad^s M(\fkm_\lz)\cong S_{\az(\lz)_s}.$$
 Then $\oln{\bfd(\lz)_s}=\az(\lz)_s$ for $1\leq s\leq \ell$. Applying the above proposition to $\fkm_\lz$ gives
the following result.

\begin{cor} \label{image-m-lambda}
 {\rm (1)} Let $\lz\in\Pi$ and keep the notation above. Then
$$\gz_{\bfd(\lz)}(\wt u_{\fkm_\lz})\in v^{\thz(\lz)}\wt u_{\fkm^\infty_\lz} +
\sum_{{\frak z}\in\fkM_\infty,\,\scrF({\frak z})<_\dg \fkm_\lz}\cZ \wt u_{\frak z},$$
 where $\thz(\lz)=\sum_{s<t}\kappa(\bfd(\lz)_s, \bfd(\lz)_t)-\sum_{s=1}^\ell h(\bfd(\lz)_s)$.

{\rm (2)} Let $\bfd\in\bbN I_\infty$ with $\oln\bfd=\az(\lz)$. If $\bfd=\tau^{rn}(\bfd(\lz))$ for some $r\in\bbZ$, then
$$\gz_{\bfd}(\wt u_{\fkm_\lz})\in v^{\thz(\lz)}\wt u_{\tau^{rm}(\fkm^\infty_\lz)} +
\sum_{{\frak z}\in\fkM_\infty,\,\scrF({\frak z})<_\dg \fkm_\lz}\cZ \wt u_{\frak z}.$$
Otherwise,
 $$\gz_{\bfd}(\wt u_{\fkm_\lz})\in \sum_{{\frak z}\in\fkM^\bfd_\infty,\,\scrF({\frak z})<_\dg \fkm_\lz}\cZ \wt u_{\frak z}.$$
\end{cor}

 In the following we briefly recall the canonical basis of $\cH(\dt)$ for $\dt=\dt_n$ or $\dt_\infty$.
By \cite{L90} and \cite[Prop.~7.5]{VV}, there is a semilinear ring
involution $\iota:\cH(\dt)\to \cH(\dt)$ taking $v\mapsto v^{-1}$ and
$\wt u_\bfd\mapsto \wt u_\bfd$ for all $\bfd\in\bbZ I_n$. It is often called the bar-involution,
usually written as $\bar x=\iota(x)$. The canonical basis
(or the global crystal basis in the sense of Kashiwara)
${\bf B}:=\{b_\fkm \mid\fkm\in\fkM\}$ for $\cH(\dt)$ (at $v=\infty$) can be
characterized as follows:
\begin{equation} \label{can-basis-elt}
\overline{b}_\fkm=b_\fkm,\;\; b_\fkm\in \wt u_\fkm +\sum_{\fkp<_\dg \fkm}v^{-1}\bbZ[v^{-1}]\wt u_\fkp;
\end{equation}
see \cite{L90}. The canonical basis elements $b_\fkm$ also admit a geometric characterization
given in \cite{L91,VV}. Let $H^i_{\cO_\fkp}(IC_{\cO_\fkm})$
be the stalk at a point of $\cO_\fkm$ of the $i$-th intersection
cohomology sheaf of the closure $\overline {\cO_\fkp}$ of
$\cO_\fkp$. Then
$$b_\fkm=\sum_{i\in\bbN\atop \fkp\leq_\dg\fkm}v^{i-\text{dim}\,\cO_\fkm+\text{dim}\,\cO_\fkp}
\dim H_{\cO_\fkp}^i(IC_{\cO_\fkm})\wt u_\fkp.$$
 For the cyclic quiver case, by \cite{L92}, the subset of $\bf B$
$${\bf B}^{\rm ap}:=\{b_\fkm\mid \fkm\in \fkM_n^{\rm ap}\}$$
 is the canonical basis of $\cC(\dt_n)$, where $\fkM_n^{\rm ap}$ denotes the set of aperiodic
multisegments, that is, those multisegments $\fkm=\sum_{i\in I_n,\,l\geq 1} m_{i,l}[i,l)$ satisfying that
for each $l\geq 1$, there is some $i\in I_n$ such that $m_{i,l}=0$. In other words, ${\bf B}^{\rm ap}$
is the canonical basis of $\bfU_v^\pm(\widehat{\frak{sl}}_n)$. Note that for each
$\lz=(\lz_1,\ldots,\lz_m)\in\Pi$, the corresponding multisegment $\fkm_\lz$ is aperiodic if and
only if $\lz$ is $n$-regular which, by definition, satisfies $\lz_s>\lz_{s+n-1}$ for $1\leq s\leq s+n-1\leq m$.

\section{Double Ringel--Hall algebras and highest weight modules}

\def\dHallp{{\scrD(\dt)^{\geq 0}}}
\def\dHalln{{\scrD(\dt)^{\leq 0}}}
\def\bfd{{\bf d}}
\def\fka{{\mathfrak a}}

In this section we follows \cite{X97, DDF} to define the double Ringel--Hall algebra $\scrD(\dt)$
of the quiver $\dt=\dt_n$ or $\dt_\infty$ and study the irreducible highest weight modules
of $\scrD(\dt_n)$ associated with integral dominant weights in terms of a quantized generalized Kac--Moody
algebra.

The Ringel--Hall algebra $\bfcH(\dt)$ of $\dt$ can be extended to a Hopf algebra $\dHallp$ which is a $\bbQ(v)$-vector
space with a basis $\{u_\fkm^+K_{\az}\mid \az\in\bbZ I, \fkm\in \fkM\}$; see \cite{R90a,Gr95,X97}
or \cite[Prop.~1.5.3]{DDF}.
Its algebra structure is given by
\begin{equation}\label{mult-formula-positive}
\begin{split}
&K_{\az} K_{\bz}=K_{\az+\bz},\;\;K_\az u_\fkm^+=v^{(\bfd(\fkm),\az)}u_\fkm^+K_\az, \\
&
u_\fkm^+u_{\fkm'}^+=\sum_{\fkp\in\fkM}v^{\lan \bfd(\fkm),\bfd(\fkm')\ran}\vphi_{\fkm,\fkm'}^\fkp(v^{2}) u_\fkp^+,
\end{split}
\end{equation}
where $\fkm,\fkm'\in\fkM$ and $\az,\bz\in\bbZ I$, and its coalgebra structure is given by
\begin{equation}\label{comult-formula-positive}
\begin{split}
&\dt(u_\fkm^+)=\sum_{\fkm',\fkm''\in\fkM}v^{\lan \bfd(\fkm'),\bfd(\fkm'')\ran}
\frac{\fka_{\fkm'}(v^2)\fka_{\fkm''}(v^2)}{\fka_\fkm(v^2)}
\vphi_{\fkm',\fkm''}^\fkm(v^{2}) u_{\fkm''}^+\otimes  u_{\fkm'}^+ K_{\bfd(\fkm'')},\\
&\dt(K_\az)=K_\az\otimes K_\az,\;\;\vez(u_\fkm^+)=0\;(\fkm\neq0),\;\;\,\vez(K_\az)=1,
\end{split}
\end{equation}
where $\fkm\in\fkM$ and $\az\in\bbZ I$. We refer to \cite{X97} or \cite{DDF} for the definition of the antipode.

Dually, there is a Hopf algebra $\dHalln$ with basis
$\{K_\az u_\fkm^- \mid \az\in\bbZ I, \fkm\in \fkM\}$. In particular, the multiplication is given by
\begin{equation}\label{mult-formula-negative}
\begin{split}
&K_{\az} K_{\bz}=K_{\az+\bz},\;\;K_\az u_\fkm^-=v^{-(\bfd(\fkm),\az)}u_\fkm^-K_\az,\\
&u_\fkm^-u_{\fkm'}^-=\sum_{\fkp\in\fkM}
v^{\lan \bfd(\fkm'),\bfd(\fkm)\ran}\vphi_{\fkm',\fkm}^\fkp(v^{2}) u_\fkp^-,
\end{split}
\end{equation}
where $\fkm,\fkm'\in\fkM$ and  $\az,\bz\in\bbZ I$. The comultiplication and the counit are given by
\begin{equation}\label{comult-formula-negative}
\begin{split}
&\dt(u_\fkm^-)=\sum_{\fkm',\fkm''\in\fkM}v^{\lan \bfd(\fkm'),\bfd(\fkm'')\ran}
\frac{\fka_{\fkm'}\fka_{\fkm''}}{\fka_\fkm}
\vphi_{\fkm',\fkm''}^\fkm(v^{2}) u_{\fkm''}^-K_{-\bfd(\fkm')}\otimes u_{\fkm'}^-,\\
&\dt(K_\az)=K_\az\otimes K_\az,\;\;\vez(u_\fkm^-)=0\;(\fkm\neq0),\;\;\,\vez(K_\az)=1,
\end{split}
\end{equation}
where $\az\in\bbZ I$ and $\fkm\in\fkM$.

It is routine to check that the bilinear form
$\psi:\dHallp\times \dHalln\ra \bbQ(v)$ defined by
\begin{equation}\label{skew-pairing}
\psi(K_{\az}u_\fkm^+ ,K_{\bz}u_{\fkm'}^-)=v^{(\az,\bz)-\lr{\bfd(\fkm),\bfd(\fkm)}+2d(\fkm)}
\frac{\dz_{\fkm,\fkm'}}{\fka_\fkm(v^{2})}
\end{equation}
 is a skew-Hopf pairing in the sense of \cite{Jo95}; see, for example, \cite[Prop.~2.1.3]{DDF}.

Following \cite{X97} or \cite[\S2.1]{DDF}, with the triple $(\dHallp,\dHalln,\psi)$
we obtain the associated reduced {\it double Ringel--Hall algebra} $\scrD(\dt)$
which inherits a Hopf algebra structure from those of $\dHallp$ and
$\dHalln$. In particular, for all elements $x\in \dHallp$ and $y\in\dHalln$,
we have  in $\dHalln$ the following relations
\begin{equation}\label{comm-rel-double-RH}
\sum \psi(x_1,y_1)  y_2 x_2 =\sum  \psi(x_2,y_2) x_1 y_1,
\end{equation}
where $\Delta(x)=\sum x_1\otimes x_2$ and $\Delta(y)=\sum y_1\otimes y_2$. Moreover, $\scrD(\dt)$ admits a triangular decomposition
\begin{equation}\label{tri-decomposition}
\scrD(\dt)=\scrD(\dt)^+\otimes \scrD(\dt)^0\otimes \scrD(\dt)^-,
\end{equation}
 where $\scrD(\dt)^\pm$ are subalgebras generated by $u_\fkm^\pm$ ($\fkm\in\fkM$), and
$\scrD(\dt)^0$ is generated by $K_\az$ ($\az\in\bbZ I$). Thus,
$\scrD(\dt)^0$ is identified with the Laurent polynomial ring $\bbQ(v)[K_i^{\pm1}: i\in I]$,
$$\aligned
\bfcH(\dt)=\cH(\dt)\otimes_\cZ\bbQ(v)&\stackrel{\sim}{\lra}\scrD(\dt)^+,\;\; u_\fkm\lmto u_\fkm^+,\\
\bfcH(\dt)^{\rm op}=\cH(\dt)^{\rm op}\otimes_\cZ\bbQ(v)&\stackrel{\sim}{\lra}\scrD(\dt)^-,\;\; u_\fkm\lmto u_\fkm^-.
\endaligned$$
  The canonical basis of $\cH(\dt)$ given in \eqref{can-basis-elt} gives the canonical bases
${\bf B}^\pm:=\{b_\fkm^\pm \mid \fkm\in\fkM\}$ of $\scrD(\dt)^\pm$ satisfying
\begin{equation} \label{can-basis-elt-pm}
b^\pm_\fkm\in \wt u_\fkm^\pm +\sum_{\fkp<_\dg \fkm}v^{-1}\bbZ[v^{-1}]\wt u^\pm_\fkp.
\end{equation}

 For $i\in I$, $\az\in\bbN I$ and $\fkm\in\fkM$, we write
$$u_i^\pm=u_{[S_i]}^\pm,\,\;u^\pm_\az=u^\pm_{[S_\az]},\,\;\text{ and }\;
\wt u_\fkm^\pm=v^{{\rm dim\,End}_\dt(M(\fkm))-{\rm dim} M(\fkm)}u_\fkm^\pm.$$
 It is known that $\scrD(\dt_\infty)$ is generated by $u_i^\pm, K_i^{\pm1}$ ($i\in \bbZ$) and is isomorphic to
$\bfU_v(\frak{sl}_\infty)$. By \cite{R93}, the $\bbQ(v)$-subalgebra of $\scrD(\dt_n)$ generated
by $u_i^\pm, K_i^{\pm1}$ ($i\in I_n=\bbZ/n\bbZ$) is isomorphic to
$\bfU'_v(\widehat{\frak{sl}}_n)$, while $\scrD(\dt_n)$ is isomorphic to $\bfU_v(\widehat{\frak{gl}}_n)$; see
\cite{Sch04,Hub2,DDF}. From now on, we write for notational simplicity,
$$\scrD(\infty)=\scrD(\dt_\infty)\;\text{ and }\;\scrD(n)=\scrD(\dt_n).$$

\begin{rems} \label{extended-double-Hall-alg}
(1) The construction of $\scrD(n)$ is slightly different from that in \cite[\S2.1]{DDF}.
In particular, the $K_i$ here play a role as $\wt K_i=K_iK_{i+1}^{-1}$ there. In particular,
they do not satisfy the equality $K_0K_1\cdots K_{n-1}=1$.

(2) We can extend $\scrD(n)$ to the $\bbQ(v)$-algebra $\widehat\scrD(n)$ by adding new generators
$D^{\pm1}$ with relations
$$DD^{-1}=1=D^{-1}D,\,K_i D=DK_i,\,DE_i=v^{\dz_{0,i}}E_i D,\,D u_\fkm^\pm=v^{-a_0}u_\fkm^\pm D$$
for all $i\in I_n$ and $\fkm\in\fkM$, where $\bfd(\fkm)=(a_i)_{i\in I_n}$. Then $\bfU_v(\widehat{\frak{sl}}_n)$
clearly becomes a subalgebra of $\widehat\scrD(n)$.
\end{rems}

As in \eqref{central elements}, define for each $t\geq 1$,
 $$\bdc_t^\pm=(-1)^tv^{-2tn}\sum_{\fkm}(-1)^{{\rm dim}\,{\rm End}(M(\fkm))}\fka_\fkm(v^{2}) u^\pm_\fkm \in\scrD(n)^\pm,$$
By Theorem \ref{Schiffmann-Hubery}, the elements
$\bdc^+_t$ and $\bdc^-_t$ are central in $\scrD(n)^+$ and $\scrD(n)^-$,
respectively. Following \cite[Sect.~4]{Hub2},
define recursively for $t\geq 1$,
$$\bdx_t^\pm=t\bdc_t^\pm-\sum_{s=1}^{t-1} \bdx_s^\pm \bdc_{t-s}^\pm\in\scrD(n)^\pm.$$
Clearly, $\bdx_t^+$ and $\bdx_t^-$ are again central elements in $\scrD(n)^+$ and $\scrD(n)^-$, respectively.
By applying \cite[Cor.~10 \& 12]{Hub1}, the $\bdx_t^\pm$ are primitive, i.e.,
$$\Delta(\bdx_t^+)=\bdx_t^+\otimes K_{t\dz}+1\otimes \bdx_t^+\;\text{ and }\;
\Delta(\bdx_t^-)=\bdx_t^-\otimes 1+K_{-t\dz}\otimes \bdx_t^-,$$
and they satisfy
$$\psi(\bdx_t^+,\bdx_s^-)=v^{2tn}\{\bdx_t,\bdx_s\}=\dz_{t,s}tv^{2tn}v^{-2tn}(1-v^{-2tn})=\dz_{t,s}t(1-v^{-2tn}).$$
 Finally, as in \cite[\S~2.2]{DDF}, we scale the elements $\bdx_t^\pm$ by setting
$$\bdz_t^\pm=\frac{v^{tn}}{v^t-v^{-t}}\bdx_t^\pm\in\scrD(n)^\pm\;\text{ for $t\geq 1$}.$$
 Then
\begin{equation}\label{primitive-z-pm}
\Delta(\bdz_t^+)=\bdz_t^+\otimes K_{t\dz}+1\otimes \bdz_t^+,\; \Delta(\bdz_t^-)=\bdz_t^-\otimes 1+K_{-t\dz}\otimes \bdz_t^-,
\end{equation}
  and
$$\psi(\bdz_t^+,\bdz_s^-)=\dz_{t,s}\frac{t(v^{2tn}-1)}{(v^t-v^{-t})^2}.$$

\begin{lem} \label{commutator-rel} {\rm (1)} For each $i\in I_n$,
$$[u_i^+,u_i^-]=\frac{K_i-K_i^{-1}}{v-v^{-1}}.$$

{\rm (2)} For $\az\in \bbN I_n$ and $t,s\geq 1$, $K_\az \bdz_t^\pm=\bdz_t^\pm K_\az$ and
\begin{equation}\label{comm-rel-central-elt}
[\bdz_t^+,\bdz_s^-]=\dz_{t,s}\, \frac{t(v^{2tn}-1)}{(v^t-v^{-t})^2}(K_{t\dz}-K_{-t\dz}).
\end{equation}
Moreover, for each $i\in I_n$ and $t\geq 1$,
$$[u_i^+,\bdz_t^-]=0=[u_i^-,\bdz_t^+].$$
\end{lem}

\begin{pf} We only prove the formula \eqref{comm-rel-central-elt}. The remaining ones are obvious.
Since $\Delta(\bdz_t^+)=\bdz_t^+\otimes K_{t\dz}+1\otimes \bdz_t^+$ and
$\Delta(\bdz_s^-)=\bdz_s^-\otimes 1+K_{-s\dz}\otimes \bdz_s^-$, we have by \eqref{comm-rel-double-RH} that
$$\aligned
&K_{t\dz}\psi(\bdz_t^+, \bdz_s^-)+ \bdz_t^+\psi(1, \bdz_s^-)+ \bdz_s^-K_{t\dz}\psi(\bdz_t^+,K_{-s\dz})+\bdz_s^-\bdz_t^+\psi(1,K_{-s\dz})\\
=& \bdz_t^+\bdz_s^-\psi(K_{t\dz},1)+\bdz_s^-\psi(\bdz_t^+,1)+\bdz_t^+K_{-s\dz}\psi(K_{t\dz}, \bdz_s^-)+K_{-s\dz}\psi(\bdz_t^+,\bdz_s^-).
\endaligned$$
 This implies that
$$[\bdz_t^+,\bdz_s^-]=\psi(\bdz_t^+,\bdz_s^-)(K_{t\dz}-K_{-s\dz})=\dz_{t,s}\frac{t(v^{2tn}-1)}{(v^t-v^{-t})^2}(K_{t\dz}-K_{-t\dz})$$
 since $\psi(1, \bdz_s^-)\!=\psi(\bdz_t^+,K_{s\dz})=\!\psi(\bdz_t^+,1)=\!\psi(K_{t\dz},\bdz_s^-)=\!0$ and
 $\psi(1,K_{s\dz})=\!\psi(K_{-t\dz},1)=\!1$.
\end{pf}

Using arguments similar to those in the proof of \cite[Th.~2.3.1]{DDF}, we obtain a presentation of $\scrD(n)$.
More precisely, $\scrD(n)$ is the $\bbQ(v)$-algebra generated by $K_i^{\pm 1}$, $u_i^+=E_i$, $u_i^-=F_i$, and
$\bdz^\pm_t$ for $i\in I_n$ and $t\geq 1$ with defining relations:

\begin{itemize}
\item[(DH1)] $K_{i}K_{j}=K_{j}K_{i},\ K_{i}K_{i}^{-1}=1=K_{i}^{-1}K_{i}$;

\item[(DH2)] $K_{i}E_j=v^{a_{ij}}E_jK_{i}$,
$K_{i}F_j=v^{-a_{ij}} F_jK_i$, $K_i \bdz^\pm_t=\bdz^\pm_t K_i$;

\item[(DH3)] $[E_i,F_j]=\delta_{i,j}\frac {K_{i}-K_{i}^{-1}}{v-v^{-1}}$, $[E_i,\bdz_t^-]=0$, $[\bdz_t^+,F_i]=0$,\\
$[\bdz_t^+,\bdz_s^-]=\dz_{t,s}\, \frac{t(v^{2tn}-1)}{(v^t-v^{-t})^2}(K_{t\dz}-K_{-t\dz})$;

\item[(DH4)]
$\displaystyle\sum_{a+b=1-c_{i,j}}(-1)^a\leb{1-c_{i,j}\atop a}\rib
E_i^{a}E_jE_i^{b}=0$ for $i\not=j$, \\ $\bdz^+_t \bdz^+_s=\bdz^+_s \bdz^+_t$, $E_i \bdz^+_t=\bdz^+_t E_i$;

\item[(DH5)] $\displaystyle\sum_{a+b=1-c_{i,j}}(-1)^a\leb{1-c_{i,j}\atop a}\rib
F_i^{a}F_jF_i^{b}=0$ for $i\not=j$,\\
 $\bdz^-_t \bdz^-_s=\bdz^-_s \bdz^-_t$, $F_i \bdz^-_t=\bdz^-_t F_i$,

\end{itemize}
 where $i,j\in I_n$ and $t,s\geq 1$.

In the following we simply identify $I_n=\bbZ/n\bbZ$ with the subset $\{0,1,\ldots,n-1\}$ of $\bbZ$.
 Let $P^\vee=(\oplus_{i\in I_n}\bbZ h_i)\oplus \bbZ d$ be
the free abelian group with basis $\{h_i\mid i\in I_n\}\cup\{d\}$. Set ${\frak h}=P^\vee\otimes_\bbZ\bbQ$ and define
$$P=\{\llz\in{\frak h}^*=\Hom_\bbQ({\frak h},\bbQ)\mid \llz(P^\vee)\subset \bbZ\}.$$
 Then $P=(\oplus_{i\in I_n}\bbZ \llz_i)\oplus \bbZ \oz$, where $\{\llz_i\mid i\in I_n\}\cup\{\oz\}$ is the
dual basis of $\{h_i\mid i\in I_n\}\cup\{d\}$. This gives rise to the Cartan datum
$(P^\vee,P,\Pi^\vee,\Pi)$ associated with the Cartan matrix $C_n=(a_{ij})$, where
$\Pi^\vee=\{h_i\mid i\in I_n\}$ is set of simple coroots and
$\Pi=\{\az_i\mid i\in I_n\}$ is the set of simple roots defined by
$$\az_i(h_j)=a_{ji},\;\az_i(d)=\dz_{0,i}\;\text{ for all $i,j\in I_n$}.$$
 Finally, let
$$P^+=\{\llz\in P\mid \llz(h_i)\geq 0,\;\forall\, i\in I_n\}=\big(\bigoplus_{i\in I_n}\bbN\llz_i\big)\oplus\bbZ\oz$$
 denote the set of dominant weights.

For each $\llz\in X$, consider the left ideal $J_\llz$ of $\scrD(n)$ defined by
$$\aligned
J_\llz&=\sum_{\fkm\in \fkM_n\backslash\{0\}}\scrD(n)u_\fkm^++\sum_{\az\in \bbZ I_n}\scrD(n)(K_\az-v^{\llz(\az)})\\
&=\sum_{\fkm\in \fkM_n\backslash\{0\}}\scrD(n)u_\fkm^++\sum_{i\in I_n}\scrD(n)(K_i-v^{\llz(h_i)}),
\endaligned$$
 where $\llz(\az)=\sum_{i\in I_n}a_i\llz(h_i)$ if $\az=\sum_{i\in I_n}a_i\vez_i\in\bbZ I_n$.
The quotient module
$$M(\llz):=\scrD(n)/J_\llz$$
 is called the Verma module which is a highest weight
module with highest vector $\eta_\llz:=1+J_\llz$. Applying the triangular decomposition
\eqref{tri-decomposition} shows that
$$\scrD(n)^-\lra M(\llz),\;x^-\lmto x^-+J_\llz$$
 is an isomorphism of $\bbQ(v)$-vector spaces. Via this isomorphism, $\scrD(n)^-$ becomes
a $\scrD(n)$-module. It is clear that $M(\llz)$ contains a unique maximal submodule $M'$. This gives
an irreducible $\scrD(n)$-module $L(\llz)=M(\llz)/M'$.

\begin{rem} By the construction, if $\llz,\llz'\in P^+$ satisfy $\llz-\llz'\in\bbZ\oz$, then
$L(\llz)=L(\llz')$. Therefore, it might be more appropriate to work with the algebra
$\widehat\scrD(n)$ defined in Remark \ref{extended-double-Hall-alg}(2).
\end{rem}

\begin{prop} \label{irred-module-rel} Let $\llz=\sum_{i\in I_n}a_i\llz_i+b \oz\in P^+$ be a
dominant weight with $\sum_{i\in I_n}a_i>0$. Then
$$L(\llz)\cong \scrD(n)^-/\big(\sum_{i\in I_n}\scrD(n)^- (u_i^-)^{a_i+1}\big).$$
\end{prop}

\begin{pf} As in \cite[Sect.~3]{DDX}, we extend the Cartan matrix $C=(a_{ij})_{i,j\in I_n}$ to
a Borcherds--Cartan matrix $\wt C=(\wt a_{ij})_{i,j\in\bbN}$ by setting
$\wt a_{ij}=a_{ij}$ for $0\leq i,j<n$ and $\wt a_{ij}=0$ otherwise. Consider the free abelian group
$\wt P^\vee=(\oplus_{i\in\bbN}\bbZ h_i)\oplus(\oplus_{i\in\bbN}\bbZ d_i)$ and define
$$\wt P=\{\thz\in (\wt P^\vee\otimes\bbQ)^\ast \mid \thz(\wt P^\vee)\subset\bbZ  \}.$$
 We then obtain a Cartan datum of type $\wt C$
 $$(\wt P^\vee, \wt P, \wt \Pi^\vee=\{h_i\mid i\in\bbN\},\wt \Pi=\{\wt\az_i\mid i\in\bbN\})$$
 where the $\wt\az_i$ are defined by
$$\wt\az_i(h_j)=\wt a_{ji}\;\text{ and }\;\wt\az_i(d_j)=\dz_{i,j},\;\;\forall\, i,j\in\bbN.$$
 Following \cite[Def.~2.1]{Kang} or \cite[Def.~1.3]{JKK}, with the above Cartan datum we have the associated quantum
 generalized Kac--Moody algebra $\bfU_v(\wt C)$ which is by definition a $\bbQ(v)$-algebra generated by
$K_i^{\pm1}, D_i^{\pm 1}, E_i, F_i$ for $i\in\bbN$ with relations; see \cite[(1.4)]{JKK} for the details.
Clearly, the subalgebra of $\bfU_v(\wt C)$
generated by $K_i^{\pm1}, D_0^{\pm1}, E_i,F_i$ for $0\leq i<n$ is isomorphic to
$\bfU_v(\widehat{\frak{sl}}_n)$.

In order to make a comparison with $\scrD(n)$, we consider the subalgebra $\wt \bfU$
of $\bfU_v(\wt C)$ generated by $K_i^{\pm1}, E_i, F_i$ for $i\in\bbN$. Then $\wt\bfU$ admits a triangular decomposition
$$\wt\bfU=\wt\bfU^-\otimes \wt\bfU^0\otimes \wt\bfU^+,$$
 where $\wt\bfU^-$, $\wt\bfU^+$, and $\wt\bfU^0$ are subalgebras generated by $F_i$, $E_i$, and $K_i^{\pm1}$
for $i\in\bbN$, respectively. In particular, $\wt\bfU^0=\bbQ(v)[K_i^{\pm1}:i\in\bbN]$.
It follows from the definition that there is a surjective algebra homomorphism $\Psi:\wt\bfU\ra\scrD(n)$ given by
$$\Psi(E_i)=\begin{cases} u_i^+, &\text{if $0\leq i<n$};\\
                        y_{i-n+1} z_{i-n+1}^+, &\text{if $i\geq n$},\end{cases}  \quad
\Psi(F_i)=\begin{cases} u_i^-, &\text{if $0\leq i<n$};\\
                        z_{i-n+1}^-, &\text{if $i\geq n$}\end{cases},   \quad \text{ and}$$
$$\Psi(K_i^{\pm1})=\begin{cases} K_i^{\pm1}, &\text{if $0\leq i<n$};\\
  K_{(i-n+1)\dz}^{\pm1}, &\text{if $i\geq n$},\end{cases} $$
where $y_t=t(v^{2tn}-1)(v-v^{-1})/(v^t-v^{-t})^2$ for $t\geq 1$; see \eqref{comm-rel-central-elt}.
Hence, each $\scrD(n)$-module can be viewed as a $\wt\bfU$-module via the homomorphism $\Psi$.
In what follows, we will identify $\wt\bfU^\pm$ with $\scrD(n)^\pm$ via $\Psi$.

As defined in \cite[Sect.~2.1]{JKK}, for each $\thz\in \wt P$, there is an associated irreducible
$\wt\bfU$-module $L(\thz)$. By \cite[Prop.~3.3]{JKK}, $L(\thz)$ is integrable if and only if
$\thz$ is dominant, that is,
$$\thz\in \wt P^+=\{\rho\in (\wt P^\vee\otimes\bbQ)^\ast \mid \rho(\wt P^\vee)\subset\bbN\}.$$
Moreover, by \cite[Cor.~4.7]{Kang},
$$L(\thz)\cong \wt\bfU^-/\big(\sum_{i\in I_n}\wt\bfU^- F_i^{\thz(h_i)+1}+\sum_{i\geq n, \thz(h_i)=0}\wt\bfU^- F_i\big).$$

Viewing the irreducible $\scrD(n)$-module $L(\llz)$ as a $\wt\bfU$-module, it is
then isomorphic to $L(\wt\llz)$, where $\wt\llz\in \wt P$ is defined by
$$\wt\llz(h_i)=\begin{cases} \llz(h_i)=a_i, &\text{if $0\leq i<n$};\\
                             (i-n+1)\sum_{0\leq j<n} a_j, &\text{if $i\geq n$}
                             \end{cases}
\;\;\text{ and }\;\;  \wt\llz(d_i)=\dz_{i,0}b.$$
 From the assumption $\sum_{i\in I}a_i>0$ it follows that $\wt\llz(h_i)>0$ for all $i\geq n$. Consequently,
$$L(\llz)\cong L(\wt\llz)\cong \wt\bfU^-/\big(\sum_{i\in I_n}\wt\bfU^- F_i^{a_i+1}\big)
= \scrD(n)^-/\big(\sum_{i\in I_n}\scrD(n)^- (u_i^-)^{a_i+1}\big).$$
\end{pf}

For each $\llz\in P$, let $L_0(\llz)$ denote the irreducible $\bfU'_v(\widehat{\frak{sl}}_n)$-module
of highest weight $\llz$. Applying Theorem \ref{Schiffmann-Hubery} gives the following result.

\begin{cor} \label{decomp-hwm} Let $\llz=\sum_{i\in I_n}a_i\llz_i+b \oz\in P^+$ with $\sum_{i\in I_n}a_i>0$.
Then $L_0(\llz)$ is the $\bfU'_v(\widehat{\frak{sl}}_n)$-submodule of $L(\llz)$ generated by the highest
weight vector $\eta_\llz$ and there is a vector space decomposition
$$L(\llz)= L_0(\llz)\otimes \bbQ(v)[\bdz_1^-,\bdz_2^-,\ldots].$$
 In particular, if $L(\llz)|_{\bfU'_v(\widehat{\frak{sl}}_n)}$ denotes the
$\bfU'_v(\widehat{\frak{sl}}_n)$-module via restriction, then
\begin{equation}\label{decomposition-L}
L(\llz)|_{\bfU'_v(\widehat{\frak{sl}}_n)}\cong \bigoplus_{m\geq 0} L_0(\llz-m\dz^*)^{\oplus p(m)},
\end{equation}
 where $\dz^*=\sum_{i\in I_n}\az_i$ and $p(m)$ is the number of partitions of $m$.
\end{cor}

\begin{pf} By Theorem \ref{Schiffmann-Hubery},
$$\scrD(n)^-=\bfU^-_v(\widehat{\frak{sl}}_n)\otimes \bbQ(v)[\bdz_1^-,\bdz_2^-,\ldots].$$
 This implies that
$$L(\llz)\cong \scrD(n)^-/\big(\sum_{i\in I_n}\scrD(n)^- (u_i^-)^{a_i+1}\big)\cong
\big(\bfU^-_v(\widehat{\frak{sl}}_n)/\big(\sum_{i\in I_n}\bfU^-_v(\widehat{\frak{sl}}_n) F_i^{a_i+1}\big)\big)
\otimes \bbQ(v)[\bdz_1^-,\bdz_2^-,\ldots].$$
By \cite[Cor.~6.2.3]{L93}, $L_0(\llz)\cong \bfU^-_v(\widehat{\frak{sl}}_n)/\big(\sum_{i\in I_n}\bfU^-_v(\widehat{\frak{sl}}_n) F_i^{a_i+1}\big)$.
Hence, $L_0(\llz)$ is the $\bfU'_v(\widehat{\frak{sl}}_n)$-submodule of $L(\llz)$ generated by $\eta_\llz$
and the desired decomposition is obtained.

For each family of nonnegative integers $\{m_t\mid t\geq 1\}$ satisfying all but finitely many $m_t$ are zero,
$L_0(\llz)\otimes \prod_{t\geq 1}(\bdz_t^-)^{m_t}$ is a $\bfU'_v(\widehat{\frak{sl}}_n)$-submodule of $L(\llz)$
since $[u_i^\pm,\bdz_t^-]=0$ for all $i\in I_n$ and $t\geq 1$. It is easy to see that
$$L_0(\llz)\otimes \prod_{t\geq 1}(\bdz_t^-)^{m_t}\cong L_0(\llz-(\sum_{t\geq1} m_t)\dz^\ast).$$
We conclude that
$$L(\llz)|_{\bfU'_v(\widehat{\frak{sl}}_n)}\cong \bigoplus_{m\geq 0} L_0(\llz-m\dz^*)^{\oplus p(m)}.$$
\end{pf}

By \cite[Th.~14.4.11]{L93}, for each $\llz\in P^+$, the canonical basis $\{b^-_\fkm \mid \fkm\in\fkM_n^{\rm ap}\}$
of $\bfU^-_v(\widehat{\frak{sl}}_n)$ gives rise to the canonical basis
$$\{b^-_\fkm \eta_\llz\not=0 \mid \fkm\in\fkM_n^{\rm ap}\}$$
of $L_0(\llz)$. On the other hand, the crystal basis theory for the quantum generalized Kac-Moody algebra $\bfU(\wt C)$ has
been developed in \cite{JKK}. Since all the $F_i$ for $i\geq n$ correspond to imaginary simple roots and
are central in $\wt\bfU^-=\scrD(n)^-$, applying the construction in \cite[Sect.~6]{JKK} shows that the set
$${\bf B}':=\big\{\big(\prod_{i\geq n} F_i^{m_i}\big) b^-_\fkm \mid \fkm\in\fkM_n^{\rm ap} \text{ and all $m_i\in\bbN$ but finitely many are zero} \big\}$$
 forms the global crystal basis of $\wt\bfU^-=\scrD(n)^-$. We remark that
$\bf B'$ does not coincide with the canonical basis ${\bf B}^-$ of $\scrD(n)^-$.
%

\section{The $q$-deformed Fock space I: $\scrD(\infty)$-module}

\def\ot{\otimes}
\def\afH{{\widehat {\bf H}}}

In this section we introduce the $q$-deformed Fock space $\llz^\infty$ from \cite{Hay} and review its module
structure over $\scrD(\infty)=\bfU_v(\frak{sl}_\infty)$ defined in \cite{MM, VV}, as well as its
$\bfU_v(\widehat{\frak{sl}}_n)$-module structure. We also provide a proof of \cite[Prop.~5.1]{VV} by using
the properties of representations of $\Delta_\infty$. Throughout this section, we identify $\scrD(\infty)$
with $\bfU_v(\widehat{\frak{sl}}_\infty)$ via taking $u_i^+\mapsto E_i$, $u_i^-\mapsto F_i$ for all $i\in I_\infty=\bbZ$.

For each partition $\lz\in\Pi$, let $T(\lz)$ denote the tableau of shape $\lz$ whose box in the intersection of the
$i$-th row and the $j$-th column is equipped with $j-i$ (The box is then said to be with color $j-i$).
For example, if $\lz=(4,2,2,1)$, then $T(\lz)$ has the form
\begin{center}
\begin{pspicture}(0,0)(6,3)
\psset{xunit=1cm,yunit=1cm,linewidth=0.5pt}
\psline(1.2,2.6)(1.8,2.6)
\psline(1.2,2)(2.4,2)
\psline(1.2,1.4)(2.4,1.4)
\psline(1.2,0.8)(3.6,0.8)
\psline(1.2,0.2)(3.6,0.2)
\psline(1.2,2.6)(1.2,0.2)
\psline(1.8,2.6)(1.8,0.2)
\psline(2.4,2)(2.4,0.2)
\psline(3,0.8)(3,0.2)
\psline(3.6,0.8)(3.6,0.2)
\uput[l](1.75,0.5){$_0$}\uput[l](2.5,0.5){$_{1}$}
\uput[l](3.1,0.5){$_{2}$} \uput[l](3.7,0.5){$_{3}$}
\uput[l](1.9,1.1){$_{-1}$} \uput[l](2.5,1.1){$_{0}$}
\uput[l](1.9,1.7){$_{-2}$} \uput[l](2.5,1.7){$_{-1}$}
\uput[l](1.9,2.3){$_{-3}$}

\end{pspicture}
\end{center}
 For given $i\in\bbZ$, a removable $i$-box of $T(\lz)$ is by definition a box with the color $i$
 which can be removed in such a way that the new tableau has the form $T(\mu)$ for some $\mu\in\Pi$.
On the contrary, an indent $i$-box of $T(\lz)$ is a box with the color $i$ which can be added to $T(\lz)$.
 For $i\in\bbZ$ and $\lz\in\Pi$, define
$$ n_i(\lz)=|\{\text{indent $i$-boxes of $T(\lz)$}\}|-|\{\text{removable $i$-boxes of $T(\lz)$}\}|.$$

 Let $\bigwedge^\infty$ be the $\bbQ(v)$-vector space with basis
$\{|\lz\rangle\mid \lz\in\Pi\}$. Following \cite[4.2]{VV}, there is a left $\bfU_v(\frak{sl}_\infty)$-module
structure on $\bigwedge^\infty$ defined by
\begin{equation}\label{infty-module-str}
K_i\cdot\blz=v^{n_i(\lz)}\blz,\;\,E_i\cdot \blz=\bnu,\;\,F_i\cdot \blz=\bmu, \;\forall\,i\in\bbZ, \lz\in\Pi,
\end{equation}
 where $\mu,\nu\in\Pi$ are such that $T(\mu)-T(\lz)$ and $T(\lz)-T(\nu)$ are a box with
 color $i$. As remarked in \cite[Sect.~2]{MM}, $\bigwedge^\infty$ is isomorphic to the basic
 representation of $\bfU_v(\frak{sl}_\infty)$ with the canonical basis
$\{\blz\mid \lz\in\Pi\}$.

\begin{lem} \label{indecom-action} {\rm(1)} For $i\in \bbZ$ and $\lz,\mu\in\Pi$, if $u_i^-\cdot \bmu=\blz$,
then there is an exact sequence
$$0\lra S_i\lra M(\fkm_\lz)\lra M(\fkm_\mu)\lra 0.$$

{\rm(2)} Let $\fkm=[i,l)$ for some $i\in\bbZ$ and $l\geq 1$. Then
$\wt u^-_{\fkm}\cdot\bempty\in\cZ\blz$ if $i\leq 0$ and $i+l-1\geq 0$ and $0$ otherwise, where
$\lz=(i+l,1^{(-i)})$. In particular, if $i=0$, then $\wt u^-_{\fkm}\cdot\bempty=\blz$.
\end{lem}

\begin{pf} (1) This follows directly from the definition.

(2) We proceed induction on $l$. The statement is trivial if $l=1$.
Suppose now $l>1$. By the definition, $M(\fkm)=S_i[l]$ with $\udim M(\fkm)=\sum_{j=i}^{i+l-1}\vez_j$.
Then
$$u^-_{i+l-1}\cdots u^-_{i+1} u^-_i =v^{1-l}u^-_{\fkm}+\sum_{{\frak z}<^\infty_{\dg}\fkm}v^{1-l}u^-_{\frak z}.$$
For each $\frak z$ with ${\frak z}<^\infty_{\dg}\fkm$, $M(\frak z)$ is decomposable. Thus, we may write
$$M({\frak z})=M({\frak y})\oplus M({\frak z}_1),$$
 where ${\frak y}\in \fkM_\infty$ and ${\frak z}_1=[j,i+l-j)$ for some $i<j\leq i+l-1$. This implies that
$$u^-_{\frak y} u^-_{{\frak z}_1}=u^-_{\frak z}.$$
By the induction hypothesis,
$$u^-_{{\frak z}_1}\cdot\bempty\in\cZ\bmu\;\;\text{if $j\leq 0$ and $i+l-1\geq 0$,}$$
 and $0$ otherwise, where $\mu=(i+l,1^{(-j)})$. Let now $j\leq 0$ and $i+l-1\geq 0$
and let $k_1,\ldots,k_{j-i}$ be a permutation of $i,i+1,\ldots,j-1$. Then
$$(u^-_{k_1}u^-_{k_2}\cdots u^-_{k_{j-i}})\cdot \bmu=0$$
unless $k_1=i,k_2=i+1\ldots, k_{j-i}=j-1$, and moreover
$$(u^-_iu^-_{i+1}\cdots u^-_{j-1})\cdot \bmu=\blz.$$
 Since $u^-_{\frak y}$ is a $\cZ$-linear combination of the monomials $u^-_{k_1}u^-_{k_2}\cdots u^-_{k_{j-i}}$,
we have $\wt u^-_{\fkm}\cdot\bempty\in\cZ\blz$.

Now let $i=0$. Then $u^-_{{\frak z}_1}\cdot\bempty=0$ for each ${\frak z}_1=[j,i+l-j)$ with $0<j\leq i+l-1$.
Hence,
$$\wt u^-_{\fkm}\cdot \bempty=v^{1-l} u^-_{\fkm}\cdot \bempty=(u^-_{l-1}\cdots u^-_1 u^-_0)\cdot\bempty
+\sum_{{\frak z}<^\infty_{\dg}\fkm}u^-_{\frak z}\cdot\bempty=\blz.$$
\end{pf}

\begin{lem} \label{zero-action} Let $\fkm=\sum_{l\geq 1}m_{i,l}[i,l)\in\fkM_\infty$ and $\lz\in\Pi$.
\begin{itemize}
\item[(1)] If there is $j\in\bbZ$ such that $\sum_{l\geq 1}m_{j,l}\geq 2$, then
$\wt u^-_{\fkm}\cdot |\lz\rangle=0$. In particular, for each $i\in\bbZ$ and $t\geq 2$, $(u^-_i)^{(t)}\cdot\blz=0$.

 \item[(2)] The element $\wt u^-_{\fkm}\cdot \blz$ is a $\cZ$-linear combination of $\bmu$ with $\mu\in\Pi$.

\end{itemize}
\end{lem}

\begin{pf} (1) For each $i\in\bbZ$, we put
$$m_i=\sum_{l\geq 1}m_{i,l}\;\text{ and }\;M_i=\bigoplus_{l\geq 1}m_{i,l}S_i[l].$$
Then $M=M(\fkm)=\oplus_{i\in\bbZ}M_i$, where all but finitely many $M_i$ are zero and
$$u^-_{\fkm}=v^{-\sum_{i>j}\lan\udim M_i,\udim M_j\ran}(\cdots u^-_{[M_{-1}]} u^-_{[M_{0}]} u^-_{[M_{1}]}\cdots).$$
Suppose there is $j\in\bbZ$ with
$m=m_j\geq 2$. Then $M_j$ admits a decomposition
$$M_j=S_j[a_1]\oplus\cdots\oplus S_j[a_{m}]\;\;\text{with $a_1\geq \cdots \geq a_{m}\geq 1$.}$$
 This implies that
$$u^-_{[S_j[a_{m}]]}\cdots u^-_{[S_j[a_1]]}=v^{b_j} u^-_{[M_j]},$$
  where $b_j=\sum_{1\leq p<q\leq m}\lr{\udim S_j[m_p],\udim S_j[m_q]}$. Hence, it suffices to show that
  for each $\mu\in\Pi$,
$$u^-_{[M_j]}\cdot\bmu=v^{-b_j}(u^-_{[S_j[a_{m}]]}\cdots u^-_{[S_j[a_1]]})\cdot\bmu=0.$$
 By the definition, $u^-_{[S_j[a_1]]}\cdot\bmu$ is a $\bbQ(v)$-linear combination of $\nu$ which are
obtained from $\mu$ by adding a $(j+r)$-box for each $0\leq r<a_1$. Thus, each such $\nu$
does not admits an indent $j$-box. Thus, $u^-_{[S_j[a_1]]}\cdot\bnu=0$ and, hence, $(u^-_{[S_j[a_{m}]]}\cdots u^-_{[S_j[a_1]]})\cdot\bmu=0$.
We conclude that $\wt u^-_{\fkm}\cdot |\lz\rangle=0$.

(2) It is known that $\wt u^-_{\fkm}$ is a $\cZ$-linear combination of monomials of divided powers
$(u_i^-)^{(t)}$ for $i\in\bbZ$ and $t\geq 1$. Since by (1), $(u_i^-)^{(t)}\cdot\bmu=0$ for
all $i\in\bbZ$, $\mu\in\Pi$ and $t\geq 2$, it follows that $\wt u^-_{\fkm}\cdot\blz$ is then a $\cZ$-linear
combination of $(u_{i_1}^-\cdots u_{i_m}^-)\cdot\blz$, where $m=\dim M(\fkm)$ and $i_1,\ldots,i_m\in\bbZ$.
By the definition, $(u_{i_1}^-\cdots u_{i_m}^-)\cdot\blz$ either is zero or lies in $\Pi$. Therefore,
$\wt u^-_{\fkm}\cdot \blz$ is a $\cZ$-linear combination of $\bmu$ with $\mu\in\Pi$.
\end{pf}

\begin{prop} \label{action-module-A-infty}
{\rm (1)} For each $\fkm\in\fkM_\infty$,
$$\wt u^-_\fkm\cdot\bempty\in \cZ\blz\;\text{ for some $\lz\in\Pi$ with $\fkm_\lz\leq^\infty_\dg \fkm$.}$$

{\rm (2)} For each $\lz\in\Pi$,
$$\wt u^-_{\fkm_\lz}\cdot\bempty=\blz.$$
 In particular, $b^-_{\fkm_\lz}\cdot\bempty=\blz.$
\end{prop}

\begin{pf} (1) If $\wt u^-_\fkm\cdot\bempty=0$, there is nothing to prove. Now suppose
$\wt u^-_\fkm\cdot\bempty\not=0$. By Lemma \ref{zero-action}(2), we write
$$\wt u^-_\fkm\cdot\bempty=\sum_{\lz\in\Pi}f_\lz(v)\blz,$$
 where all $f_\lz(v)\in\cZ$ but finitely many are zero. If $f_\lz(v)\not=0$, then $\udim M(\fkm_\lz)=\udim M(\fkm)$.
By Lemma \ref{order-coincidence}(1), such a $\lz\in\Pi$ is unique. Hence, we may suppose $\wt u^-_\fkm\cdot\bempty=f(v)\blz$
for some $0\not=f(v)\in\cZ$ and $\lz\in\Pi$. It remains to show that $\fkm_\lz\leq^\infty_\dg \fkm$.

Applying Lemma \ref{zero-action}(1) implies that
$$M=M(\fkm)=S_{i_1}[a_1]\oplus\cdots \oplus S_{i_t}[a_t],$$
 where $i_1<\cdots<i_t$ and $a_1,\ldots,a_t\geq 1$. Then
$$u^-_{[S_{i_1}[a_1]]}\cdots u^-_{[S_{i_t}[a_t]]}=v^a u^-_{\fkm},$$
  where $a=\sum_{1\leq p<q\leq t}\lr{\udim S_{i_q}[a_q],\udim S_{i_p}[a_p]}$.

We proceed induction on $t$ to show that $M(\fkm_\lz)\leq^\infty_\dg M=M(\fkm)$. If $t=1$, this follows from
Lemma \ref{indecom-action}(2). Let now $t>1$ and let $\mu\in\Pi$ be such that
$$(u^-_{[S_{i_2}[a_2]]}\cdots u^-_{[S_{i_t}[a_t]]})\cdot\bempty=g(v)\bmu \;\text{ for some $0\not=g(v)\in \cZ$.}$$
 Then $u^-_{[S_{i_1}[a_1]]}\cdot\bmu=v^af(v)g(v)^{-1}\blz$. By the induction hypothesis,
$$M(\fkm_\mu)\leq^\infty_\dg S_{i_2}[a_2]\oplus\cdots \oplus S_{i_t}[a_t].$$
  By writing $u^-_{[S_{i_1}[a_1]]}$ as a $\cZ$-linear combination of monomials of $u_i^-$'s and applying
Lemma \ref{indecom-action}(1), there exists $X\in\rep \dt_\infty$ satisfying $\udim X=\udim S_{i_1}[a_1]$
with an exact sequence
$$0\lra X\lra M(\fkm_\lz)\lra M(\fkm_\mu)\lra 0.$$
  Since $S_{i_1}[a_1]$ is indecomposable, it follows that $X\leq_\dg^\infty S_{i_1}[a_1]$. Therefore,
$$\aligned
M(\fkm_\lz)\leq^\infty_\dg M(\fkm_\mu)\ast X\leq^\infty_\dg &(S_{i_2}[a_2]\oplus\cdots \oplus S_{i_t}[a_t])\ast S_{i_1}[a_1]\\
=&(S_{i_2}[a_2]\oplus\cdots \oplus S_{i_t}[a_t])\oplus S_{i_1}[a_1]=M(\fkm),
\endaligned$$
that is, $\fkm_\lz\leq^\infty_\dg \fkm$.

(2) Write $\lz=(\lz_1,\ldots,\lz_t)$ with $\lz_1\geq\cdots\geq \lz_t\geq 1$. Since
$$M(\fkm_\lz)=S_0[\lz_1]\oplus S_{-1}[\lz_2]\oplus\cdots\oplus S_{1-t}[\lz_t],$$
 we have that
$$u^-_{[S_{1-t}[\lz_t]]}\cdots u^-_{[S_{-1}[\lz_2]]}u^-_{[S_0[\lz_1]]}
=v^c u^-_{\fkm_\lz},$$
 where
$$c=\sum_{1\leq r<s\leq t}\lr{\udim S_{1-r}[\lz_r],\udim S_{1-s}[\lz_s]}=\sum_{1\leq r<s\leq t}\dim\Hom_{\dt_\infty}(S_{1-r}[\lz_r], S_{1-s}[\lz_s]).$$
 By using an argument similar to that in the proof of Lemma \ref{indecom-action}(2), we obtain that
$$\aligned
v^{c} u^-_{\fkm_\lz}\cdot \bempty&= (u^-_{[S_{1-t}[\lz_t]]}\cdots u^-_{[S_{-1}[\lz_2]]}u^-_{[S_0[\lz_1]]})\cdot\bempty\\
&=v^{\lz_1-1}(u^-_{[S_{1-t}[\lz_t]]}\cdots u^-_{[S_{-1}[\lz_2]]})\cdot |(\lz_1)\rangle\\
&=c^{\lz_1+\lz_2-2}(u^-_{[S_{1-t}[\lz_t]]}\cdots u^-_{[S_{-2}[\lz_3]]})\cdot|(\lz_1,\lz_2)\rangle\\
&=v^{\lz_1+\cdots+\lz_t-t}|(\lz_1,\ldots,\lz_t)\rangle=v^{\lz_1+\cdots+\lz_t-t}\blz.
\endaligned$$
 Since
$$\dim\End_{\dt_\infty}(M(\fkm_\lz))=\sum_{1\leq r\leq s\leq t}\dim\Hom_{\dt_\infty}(S_{1-r}[\lz_r], S_{1-s}[\lz_s])=c+t$$
 and $\dim M(\fkm_\lz)=\lz_1+\cdots+\lz_t$, it follows that
$$\wt u^-_{\fkm_\lz}\cdot\bempty=v^{c+t-(\lz_1+\cdots+\lz_t)} u^-_{\fkm_\lz}\cdot\bempty=\blz.$$

By \eqref{can-basis-elt-pm},
$$b_{\fkm_\lz}^-\in \wt u^-_{\fkm_\lz}+\sum_{{\fkp}<_\dg\fkm_\lz}v^{-1}\bbZ[v^{-1}]\wt u^-_{\fkp}.$$
 Let $\fkp<_\dg\fkm_\lz$ and suppose $\wt u^-_{\fkp}\cdot\bempty\not=0$. By (1), there exists $\mu\in\Pi$ with $\fkm_\mu\leq_\dg\fkp$
such that $\wt u^-_{\fkp}\cdot\bempty=f(v)\bmu$ for some $f(v)\in\cZ$. Thus, $\fkm_\mu<_\dg\fkm_\lz$.
By Lemma \ref{order-coincidence}(1), $\mu=\lz$ since $\udim M(\fkm_\mu)=\udim M(\fkm_\lz)$.
This is a contradiction. Hence, $\wt u^-_{\fkp}\cdot\bempty=0$. We conclude that
$$b_{\fkm_\lz}^-\cdot\bempty=\wt u^-_{\fkm_\lz}\cdot \bempty=\blz.$$
\end{pf}

As a consequence of the proposition above, we obtain \cite[Prop.~5.1]{VV} as follows.

\begin{cor} The subspace $\cI$ of $\bfU^-_v(\frak{sl}_\infty)$ spanned by
$b^-_\fkm$ with $\fkm\in\fkM-\{\fkm_\lz\mid\lz\in\Pi\}$
is a left ideal of $\bfU^-_v(\frak{sl}_\infty)$. Moreover, the map
$$\bfU^-_v(\frak{sl}_\infty)/\cI\lra \mbox{$\bigwedge^\infty$},\;\,b^-_{\fkm_\lz}+\cI\lmto \blz,\;\lz\in\Pi$$
 is an isomorphism of $\bfU^-_v(\frak{sl}_\infty)$-modules.
\end{cor}

\begin{pf} This follows from Proposition \ref{action-module-A-infty}(2) and the fact that the set
$$\{b_\fkm^-\cdot\bempty\not=0\mid \fkm\in\fkM_\infty\}$$
 is a basis of $\bigwedge^\infty$; see \cite[Th.~14.4.11]{L93}.
\end{pf}

Finally, for $i\in\bbZ$ and $\lz\in\Pi$, put
$$n_i^-(\lz)=\sum_{j<i,\,j\in\bar i}n_j(\lz),\,\; n_i^+(\lz)=\sum_{j>i,\,j\in\bar i}n_j(\lz),\,\;
\text{ and }\; n_{\bar i}(\lz)=\sum_{j\in\bar i}n_j(\lz).$$
 By \cite{Hay,MM}, there is a $\bfU_v(\widehat{\frak{sl}}_n)$-module structure on $\bigwedge^\infty$ defined by
\begin{equation}\label{affine-sln-module-str}
K_{\bar i}\cdot \blz=v^{n_{\bar i}(\lz)}\blz,\;\,E_{\bar i}\cdot\blz=\sum_{j\in \bar i} v^{n_j^-(\lz)}E_j\cdot\blz,\;\,
F_{\bar i}\cdot\blz=\sum_{j\in \bar i}v^{-n_j^+(\lz)}F_j\cdot\blz,
\end{equation}
 where $\bar i\in I_n=\bbZ/n\bbZ$.

\section{The $q$-deformed Fock space II: ${\scrD}(n)$-module}

\def\bfi{{\bf i}}

\def\bfn{{\bf n}}

In this section we first recall the left ${\scrD}(n)^{\leq 0}$-module structure on the Fock space
$\bigwedge^\infty$ defined by Varagnolo and Vasserot in \cite{VV} and then extend their construction to
obtain a $\scrD(n)$-module structure on $\bigwedge^\infty$.

For each $x=\sum_{\fkm} x_\fkm u_\fkm\in\cH(\dt)$ with $\dt=\dt_n$ or $\dt_\infty$, we write
$$x^\pm=\sum_\fkm x_\fkm u_\fkm^\pm\in\scrD(\dt)^\pm.$$
Then for each $\bfd\in \bbN I_\infty$, the map $\gz_\bfd: \cH(\dt_n)_{\bar\bfd}\ra \cH(\dt_\infty)_{\bfd}$
defined in Section 3 induces $\bbQ(v)$-linear maps
$$\gz_\bfd^\pm: \scrD(n)^\pm_{\bar\bfd}\lra \scrD(\infty)^\pm_{\bfd}$$
such that $\gz_\bfd^\pm (x^\pm)=(\gz_\bfd(x))^\pm$ for each $x\in\cH(\dt_\infty)$.

Following \cite[6.2]{VV}, for each $\bar i\in I_n=\bbZ/n\bbZ$, $\lz\in\fkM_n$ and $x\in {\scrD}(n)^-_\az$, define
\begin{equation} \label{negative-action-Fock-space}
K_{\bar i}\cdot\blz=v^{n_{\bar i}(\lz)}\blz\;\text{ and }\; x\cdot \blz
=\sum_{\bfd}\big(\gz^-_\bfd(x)K_{-\bfd'}\big)\cdot\blz,
\end{equation}
 where the sum is taken over all $\bfd\in\bbN I_\infty$ such that $\bar\bfd=\az$ and
$\bfd'=\sum_{i>j,\,\bar i=\bar j}d_j\vez_i$. By \cite[Cor.~6.2]{VV}, this defines a
left ${\scrD}(n)^{\leq 0}$-module structure on $\bigwedge^\infty$ which extends the
Hayashi action of $\bfU^{\leq 0}_v({\widehat{\frak{sl}}}_n)$ on $\bigwedge^\infty$
defined in \eqref{affine-sln-module-str}.

Dually, for each $\lz\in\Pi$ and $x\in {\scrD}(n)^+_\az$, define
\begin{equation} \label{positive-action-Fock-space}
x\cdot \blz=\sum_{\bfd}\big(\gz^+_\bfd(x)K_{\bfd''}\big)\cdot \blz,
\end{equation}
 where the sum is taken over all $\bfd\in\bbN I_\infty$ such that $\bar\bfd=\az$ and
$\bfd''=\sum_{i<j,\,\bar i=\bar j}d_j\vez_i$.

\begin{prop} The formula \eqref{positive-action-Fock-space} defines a
left ${\scrD}(n)^{\geq 0}$-module structure on $\bigwedge^\infty$ which extends the
Hayashi action of $\bfU^{\geq 0}_v({\widehat{\frak{sl}}}_n)$ on $\bigwedge^\infty$.
\end{prop}

\begin{pf} Let $x\in {\scrD}(n)^+_\az$ and $y\in {\scrD}(n)^+_\bz$, where $\az,\bz\in\bbN I_n$.
By the definition, we have, on the one hand, that
$$(xy)\cdot \blz=\sum_{\bfd}\big(\gz^+_\bfd(xy)K_{\bfd''}\big)\cdot \blz$$
 and, on the other hand, that
$$x\cdot (y\cdot \blz)=\sum_{\bfa,\bfb}\big(\gz^+_\bfa(x)K_{\bfa''}\gz^+_\bfb(y)K_{\bfb''}\big)\cdot \blz,$$
 where the sum is taken over all $\bfa,\bfb\in\bbN I_\infty$ such that $\bar\bfa=\az$ and $\bar\bfb=\bz$.

Since $K_{\bfa''}\gz^+_\bfb(y)=v^{(\bfa'',\bfb)}\gz^+_\bfb(y)K_{\bfa''}$, we obtain that
$$x\cdot (y\cdot\blz)=\sum_{\bfd}\sum_{\bfa+\bfb=\bfd}v^{(\bfa'',\bfb)}\big(\gz^+_\bfa(x)\gz^+_\bfb(y)K_{\bfd''}\big)\cdot \blz.$$
 By the definition,
$$(\bfa'',\bfb)=\big(\sum_{i<j,\,\bar i=\bar j}a_j\vez_i,\sum_{i}b_i\vez_i\big)=
\sum_{i<j,\,\bar i=\bar j}b_i(2a_j-a_{j-1}-a_{j+1})=\kappa(\bfa,\bfb).$$
 Applying Lemma \ref{compatibility-gamma}(2) gives that
$$(xy)\cdot\blz=x\cdot (y\cdot\blz).$$
 Hence, $\bigwedge^\infty$ becomes a left ${\scrD}(n)^{\geq 0}$-module.

For each $\bar i\in I_n=\bbZ/n\bbZ$ and $\lz\in\Pi$, we have
$$u^+_{\bar i}\cdot\blz=\sum_{j\in \bar i}(u_j^+K_{-\vez_j''})\cdot\blz.$$
 Since $\vez_j''=\sum_{l<j,\,\bar l=\bar j}\vez_l$ for each $j\in\bar i$, it follows that
$$K_{\vez_j''}\cdot\blz=\prod_{l<j,\,\bar l=\bar j}K_{\vez_l}\cdot\blz
=v^{\sum_{l<j,\bar l=\bar j}n_l(\lz)}\blz=v^{n_j^-(\lz)}\blz.$$
 This implies that
$$u^+_{\bar i}\cdot\blz=\sum_{j\in \bar i} v^{n_j^-(\lz)}u^+_j\cdot\blz,$$
 which coincides with the formula for $E_{\bar i}\cdot\blz$ in \eqref{affine-sln-module-str},
as required.
\end{pf}

The main purpose of this section is to prove that formulas \eqref{negative-action-Fock-space} and
\eqref{positive-action-Fock-space} indeed define a $\scrD(n)$-module structure on $\bigwedge^\infty$.
The strategy is to pass to the semi-infinite $v$-wedge spaces defined in \cite{KMS}.

Let $\ooz$ denote the $\bbQ(v)$-vector space with basis $\{\oz_i\mid i\in\bbZ\}$.
By \cite[Prop.~3.5]{DDF}, $\ooz$ admits a $\scrD(n)$-module structure defined by
\begin{equation}\label{dot-action-u_A}
\aligned
u_i^+ \cdot\oz_s=\dz_{i+1,\bar s}\oz_{s-1},\;\; u_i^- \cdot\oz_s=\dz_{i,\bar s}\oz_{s+1}\\
K_i^{\pm1}\cdot \oz_s=v^{\pm\dz_{i,\bar s}\mp\dz_{i+1,\bar s}}\oz_s,\;\;\bdz_m^\pm\cdot \oz_s=\oz_{s\mp mn}
\endaligned\end{equation}
for all $i\in I_n$ and $s,m\in\bbZ$ with $m\geq 1$. In particular, $K^{\pm1}_\dz\cdot \oz_s=\oz_s$ for each $s\in\bbZ$.
This is an extension of the $\bfU'_v({\widehat{\frak{sl}}}_n)$-action on $\ooz$ defined in \cite[1.1]{KMS} as well
as an extension of the $\scrD(n)^{\leq 0}$-action on $\ooz$ defined in \cite[8.1]{VV}; see \cite[3.5]{DDF}.

For a fixed positive integer $r$, consider the $r$-fold tensor product $\ooz^{\ot r}$
which has a basis
$$\{\oz_\bfi=\oz_{i_1}\ot\cdots \ot\oz_{i_r} \mid\bfi=(i_1,\ldots,i_r)\in\bbZ^r\}.$$
 The Hopf algebra structure of $\scrD(n)$ induces a $\scrD(n)$-module
structure on the $r$-fold tensor product $\ooz^{\ot r}$. By \eqref{primitive-z-pm}, we have for each $t\geq 1$,
\begin{equation}\label{r-fold-z-pm}
\aligned
\Delta^{(r-1)}(\bdz^+_t)&=\sum_{s=0}^{r-1} \underbrace{1\otimes\cdots\otimes 1}_{s}\otimes \bdz^+_t\otimes
\underbrace{K_{t\dz}\otimes\cdots\otimes K_{t\dz}}_{r-s-1}\;\text{ and }\\
\Delta^{(r-1)}(\bdz^-_t)&=\sum_{s=0}^{r-1} \underbrace{K_{-t\dz}\otimes\cdots\otimes K_{-t\dz}}_{s}\otimes \bdz^-_t\otimes
\underbrace{1\otimes\cdots\otimes 1}_{r-s-1}.
\endaligned
\end{equation}
This implies particularly that for each $t\geq 1$ and $\oz_\bfi=\oz_{i_1}\ot\cdots\ot\oz_{i_r}\in \ooz^{\ot r}$,
\begin{equation}\label{action-central-elt}
\bdz^\pm_t\cdot\oz_\bfi=\sum_{s=1}^r \oz_{i_1}\ot\cdots\ot
\oz_{i_{s-1}}\ot\oz_{i_s\mp tn}\ot\oz_{i_{s+1}}\ot\cdots\ot\oz_{i_r}.
\end{equation}

%

By \eqref{comult-formula-positive} and \eqref{comult-formula-negative}, for each $\az\in\bbN I_n$, we have
\begin{equation} \label{comulti-ss-module}  \begin{split}
 &\Delta^{(r-1)}(\wt u_\az^+)
 =\sum_{\az=\az^{(1)}+\cdots+\az^{(r)}}
 v^{\sum_{s>t}\lan\az^{(s)},\az^{(t)}\ran}\times\\
 &\hspace{4cm}\wt u_{\az^{(1)}}^+\otimes\wt u_{\az^{(2)}}^+K_{\az^{(1)}}\otimes\cdots
 \otimes\wt u_{\az^{(r)}}^+K_{(\az^{(1)}+\az^{(2)}+\cdots+\az^{(r-1)})},\\
 &\Delta^{(r-1)}(\wt u_\az^-)
 =\sum_{\az=\az^{(1)}+\cdots+\az^{(r)}}
 v^{\sum_{s>t}\lan\az^{(s)},\az^{(t)}\ran}\times\\
 &\hspace{4cm}\wt u_{\az^{(1)}}^- K_{-(\az^{(2)}+\cdots+\az^{(r)})}
 \otimes \cdots \otimes u_{\az^{(r-1)}}^-K_{-\az^{(r)}}
 \otimes\wt u_{\az^{(r)}}^-.
\end{split}
\end{equation}
 This gives the the following lemma; see \cite[Lem.~8.3]{VV} and \cite[Cor.~3.5.8]{DDF}.

\begin{lem} Let $\az\in\bbN I_n$ and $\bfi=(i_1,\ldots,i_r)\in\bbZ^r$. Then
\begin{equation} \label{formula-positive-action}
\tu_\az^+\cdot\oz_\bfi=\sum_{\bfn}v^{c^+(\bfi,\bfi-\bfn)}\oz_{\bfi-\bfn},
\end{equation}
where the sum is taken over the sequences ${\bfn}=(n_1,\ldots,n_r)\in\{0,1\}^r$ satisfying
$\az=\sum_{s=1}^r n_s\vez_{\overline{i_s-1}}$ and
$$c^+(\bfi,\bfi-{\bfn})=\sum_{1\leq s<t\leq r}n_s(n_t-1)\lan \vez_{\bar i_t},\vez_{\bar i_s}\ran;$$

\begin{equation}\label{formula-negative-action}
\tu_\az^-\cdot\oz_\bfi=\sum_{\bf n}v^{c^-(\bfi,\bfi+\bf n)}\oz_{\bfi+\bf n},
\end{equation}
 where the sum is taken over the sequences ${\bf n}=(n_1,\ldots,n_r)\in\{0,1\}^r$ satisfying
$\az=\sum_{s=1}^rn_s\vez_{\bar i_s}$ and
$$c^-(\bfi,\bfi+{\bf n})=\sum_{1\leq s<t\leq r}n_t(n_s-1)\lan \vez_{\bar i_t},\vez_{\bar i_s}\ran.$$

\end{lem}

On the other hand, let $\afH(r)$ be the Hecke algebra of affine symmetric group of type $A$ which
is by definition a $\bbQ(v)$-algebra with generators
$T_i$ and $X_j$ for $i=1,\ldots,r-1$, $j=1,\ldots,r$ and relations:
$$\aligned
 & (T_i+1)(T_i-v^{2})=0,\\
 & T_iT_{i+1}T_i=T_{i+1}T_iT_{i+1},\;\;T_iT_j=T_jT_i\;(|i-j|>1),\\
 & X_iX_i^{-1}=1=X_i^{-1}X_i,\;\; X_iX_j=X_jX_i,\\
 & T_iX_iT_i=v^{2} X_{i+1},\;\;X_jT_i=T_iX_j\;(j\not=i,i+1).
\endaligned$$
 This is the so-called {\it Bernstein presentation} of $\afH(r)$.

 By \cite[Sect.~8.2]{VV}, there is a right $\afH(r)$-module structure on $\ooz^{\ot r}$ defined by
\begin{equation}\label{afH action}
\aligned
{}& \oz_{\bf i}\cdot X_t=\oz_{i_1}\cdots\oz_{i_{t-1}}\oz_{i_t-n}\oz_{i_{t+1}}\cdots\oz_{i_r},\\
& {\oz_{\bf i}\cdot T_k=\left\{\begin{array}{ll} v^{2}\oz_{\bf
i},\;\;&\text{if $i_k=i_{k+1}$;}\\
v\oz_{{\bf i} s_k},\;\;&\text{if $-n<i_k<i_{k+1}\leq 0$;}\\
v\oz_{{\bf i} s_k}+(v^{2}-1)\oz_{\bf i},\;\;&\text{if
$-n<i_{k+1}<i_k\leq 0$,}
\end{array}\right.}
\endaligned
\end{equation}
 where ${\bf i}=(i_1,\ldots, i_r)\in \bbZ^r$, $\oz_{\bf i}=\oz_{i_1}\ot\cdots\ot\oz_{i_r}$ and
$$\oz_{{\bf i}s_k}=\oz_{i_1}\ot\cdots\ot \oz_{i_{k+1}}\ot\oz_{i_k}\ot\cdots\ot\oz_{i_r}.$$
 Following \cite[Lem.~8.2]{VV} and \cite[Prop.~3.5.5]{DDF}, the tensor space $\ooz^{\ot r}$
is indeed a $\scrD(n)$-$\afH(r)$-bimodule. Set
$$\Xi^r=\sum_{i=1}^{r-1}\im(1+T_i)\subseteq \ooz^{\ot r},$$
which is clearly a $\scrD(n)$-submodule of $\ooz^{\ot r}$. Thus, the quotient space
$\ooz^{\ot r}/\Xi^r$ becomes a $\scrD(n)$-module. For each $\bfi=(i_1,\ldots,i_r)\in\bbZ^r$, write
$$\wedge\oz_\bfi=\oz_{i_1}\wedge\ldots \wedge\oz_{i_r}=\oz_\bfi+\Xi^r\in \ooz^{\ot r}/\Xi^r.$$
 By \cite[Prop.~1.3]{KMS}, the set
$$\{\wedge\oz_\bfi\mid i_1>\cdots>i_r\}$$
 forms a basis of $\ooz^{\ot r}/\Xi^r$.

For each $m\in\bbZ$, let ${\scr B}_m$ denote the set of sequences $\bfi=(i_1,i_2,\ldots)\in\bbZ^\infty$
satisfying that $i_s=m-s+1$ for $s\gg 0$, and set ${\scr B}_\infty=\cup_{m\in\bbZ}{\scr B}_m$.
As in \cite[Sect.~10.1]{VV}, let $\ooz^\infty$ denote the space spanned by semi-infinite monomials
$$\oz_\bfi=\oz_{i_1}\ot\oz_{i_2}\ot\cdots,\;\text{ where $\bfi=(i_1,i_2,\ldots)\in{\scr B}_\infty$.}$$
 Then the affine Hecke algebra $\afH(\infty)$ acts on $\ooz^\infty$
via the formulas in \eqref{afH action}. Set
$$\Xi^\infty=\sum_{i=1}^\infty \im(1+T_i)\subseteq \ooz^\infty.$$
 For each $\bfi=(i_1,i_2,\ldots)\in{\scr B}_\infty$ as above, write
$$\wedge\oz_\bfi=\oz_{i_1}\wedge\oz_{i_2}\wedge\cdots =\oz_\bfi+\Xi^\infty\in \ooz^{\infty}/\Xi^\infty.$$
 By \cite[Prop.~1.4]{KMS}, the $\bfU'_v({\widehat{\frak{sl}}}_n)$-module structure on $\ooz^{\ot r}/\Xi^r$
induces a $\bfU'_v({\widehat{\frak{sl}}}_n)$-module structure on $\ooz^{\infty}/\Xi^\infty$. Moreover,
the map
$$\kappa:\mbox{$\bigwedge^\infty$}\lra \ooz^{\infty}/\Xi^\infty,\; \blz\longmapsto \wedge\oz_{\bfi_\lz}$$
is an injective homomorphism of $\bfU'_v({\widehat{\frak{sl}}}_n)$-modules.

  Following \cite[1.4]{KMS}, for each $m\in\bbZ$, write
$$|m\ran=\oz_{m}\wedge\oz_{m-1}\wedge\oz_{m-2}\wedge\cdots.$$
 Clearly, for each $\bfi=(i_1,i_2,\ldots)\in{\scr B}_m$, there exists a sufficiently large $N$ such that
$$\oz_\bfi=(\oz_{i_1}\wedge\cdots\wedge \oz_{i_N})\wedge|m-N\ran.$$

 By \cite[Lem.~2.2]{KMS} and \eqref{comulti-ss-module}, for given $\az\in\bbN I$ and $\bfi\in{\scr B}_m$,
there is $t\gg 0$ such that
$$u_\az^-\cdot (\wedge\oz_\bfi)=\big(u_\az^-\cdot(\oz_{i_1}\wedge\cdots\wedge\oz_{i_t})\big)\wedge|m-t\ran.$$
 Hence, the $\scrD(n)^{\leq0}$-module structure on $\ooz^{\ot r}/\Xi^r$
induces a $\scrD(n)^{\leq0}$-module structure on $\ooz^{\infty}/\Xi^\infty$; see \cite[Sect.~10.1]{VV}. Moreover,
by \cite[Lem.~10.1]{VV}, the map $\kappa:\mbox{$\bigwedge^\infty$}\ra \ooz^{\infty}/\Xi^\infty$
is a $\scrD(n)^{\leq0}$-module homomorphism.

Dually, for each given $\bfi\in{\scr B}_m$, there is $t\gg 0$ such that
$$u_\az^+\cdot (\wedge\oz_\bfi)=\big(u_\az^+\cdot(\oz_{i_1}\wedge\cdots\wedge\oz_{i_t})\big)\wedge\big(K_{\az}\cdot |m-t\ran\big).$$
 Thus, $\ooz^{\infty}/\Xi^\infty$ becomes a left $\scrD(n)^{\geq 0}$-module as well. We have the following
 result.

\begin{prop} \label{Hall-algebra-action-FS} The map $\kappa$ is a $\scrD(n)^{\geq0}$-module homomorphism.
\end{prop}

\begin{pf} We need to show that for each $\lz\in\Pi$ and $\az\in\bbN I_n$,
$$\kappa(\wt u_\az^+\cdot\blz)=\tu_\az^+(\wedge\oz_{\bfi_\lz}).$$
 For simplicity, write $\bfi=\bfi_\lz$. By \eqref{positive-action-Fock-space},
$$\tu_\az^+\cdot \blz=\sum_{\bfd}(\gz^+_\bfd(\tu_\az^+)K_{\bfd''})\cdot \blz=\sum_{\bfd}v^{-h(\bfd)} (\tu_\bfd^+ K_{\bfd''})\cdot\blz,$$
 where the sum is taken over all $\bfd\in\bbN I_\infty$ such that $\bar\bfd=\az$ and
$h(\bfd)=\sum_{i<j,\bar i=\bar j}d_i(d_{j+1}-d_j)$.

For each fixed $\bfd=(d_i)\in \bbN I_\infty$ with $\bar\bfd=\az$, we have
$$\tu_\bfd^+=\cdots \tu^+_{d_1\vez_1}\tu^+_{d_0\vez_0}\tu^+_{d_{-1}\vez_{-1}}\cdots=\prod_{i\in\bbZ}\tu^+_{d_i\vez_i}.$$
 By the definition, $\tu_\bfd^+\cdot\blz\not=0$ implies that
$$\bfd=\sum_{s\geq 1}n_s\vez_{i_s-1},$$
 where $n_s\in\{0,1\}$ for all $s\geq 1$. Moreover, if this is the case, then
$$\tu_\bfd^+\cdot\blz=|\mu_\bfn\ran,$$
 where $\bfn=(n_1,n_2,\ldots)$ and $\mu_\bfn=\mu\in\Pi$ is determined by $\bfi_\mu=\bfi-\bfn$.
Therefore, for $\bfd\in\bbN I_\infty$ with $\bfd=\sum_{s\geq 1}n_s\vez_{i_s-1}$, we have that
$$ \aligned
K_{\bfd''}\cdot\blz&=\prod_{\bar i_s=\bar i_t,\,i_s>i_t}K_{i_t-1}^{n_s}\cdot\blz
=\big(\prod_{\bar i_s=\bar i_t,\,i_s>i_t}v^{n_s\sum_{l\in\bbZ}(\dz_{i_t-1,i_l}-\dz_{i_t-1,i_l-1})}\big)\cdot\blz\\
&=\prod_{i_s>i_t}v^{n_s(\dz_{\bar i_s,\overline{i_t+1}}-\dz_{\bar i_s,\bar i_t})}\blz
=\big(v^{-\sum_{i_s>i_t}n_s\lan \vez_{\bar i_t},\vez_{\bar i_s}\ran}\big)\blz
\endaligned$$
 and
$$h(\bfd)=\sum_{i_s>i_t}-n_sn_t(\dz_{\bar i_s,\bar i_t}-\dz_{\bar i_s,\overline{i_t+1}})
=-\sum_{i_s>i_t} n_sn_t\lan \vez_{\bar i_t},\vez_{\bar i_s}\ran.$$
 Since $i_s>i_t$ if and only if $s<t$, we conclude that
$$\tu_\az^+\cdot \blz=\sum_{\bfn}v^{x(\bfn)}|\mu_\bfn\ran,$$
 where the sum is taken over the sequences ${\bf n}=(n_1,n_2,\ldots)\in\{0,1\}^\infty$ satisfying
$\az=\sum_{s=1}^rn_s\vez_{\overline{i_s-1}}$ and
$$x(\bfn)=\sum_{1\leq s<t}n_s(n_t-1)\lan \vez_{\bar i_t},\vez_{\bar i_s}\ran.$$
 This together with \eqref{formula-positive-action} implies that
$$\kappa(\tu_\az^+\cdot\blz)=\sum_{\bfn}v^{x(\bfn)}\kappa(|\mu_\bfn\ran)=\sum_{\bfn}v^{x(\bfn)}\wedge\oz_{\bfi-\bfn}
=\tu_\az^+(\wedge\oz_\bfi)=\tu_\az^+(\kappa(\blz)).$$
 This finishes the proof.
\end{pf}

As a consequence of the results above, to prove that the formulas \eqref{negative-action-Fock-space} and
\eqref{positive-action-Fock-space} define a $\scrD(n)$-module structure on $\bigwedge^\infty$, it suffices
to show that the $\scrD(n)^{\leq 0}$-module and $\scrD(n)^{\geq 0}$-module structures on
$\ooz^{\infty}/\Xi^\infty$ define a $\scrD(n)$-module structure. In other words, we need to show that
the actions of $K_i^{\pm1}, u_i^+, u_i^-$ ($i\in I_n$) and $\bdz_s^+,\bdz_s^-$ ($s\geq 1$) on $\ooz^{\infty}/\Xi^\infty$ satisfy the
relations (DH1)--(DH5) in Section 4. In the following we only check the relations
$$[\bdz_t^+,\bdz_s^-]=\dz_{t,s}\, \frac{t(v^{2tn}-1)}{(v^t-v^{-t})^2}(K_{t\dz}-K_{-t\dz}).$$
 The other relations either follow from \cite{KMS} or can be checked directly.

By \cite[\S2]{KMS}, for each $t\geq 1$, there are Heisenberg operators
$B^\pm_t:\ooz^{\infty}/\Xi^\infty\ra \ooz^{\infty}/\Xi^\infty$ taking
$$B_t(\wedge\oz_{\bfi})\lmto \sum_{s=1}^\infty \wedge\oz_{\bfi\mp tn{\bf e}_s},$$
 where $\bfi\in{\scr B}_\infty$ and ${\bf e}_s=(\dz_{i,s})_{i\geq 1}\in\bbZ^\infty$.
Note that for each $\bfi\in{\scr B}_\infty$, $\wedge\oz_{\bfi\mp tn{\bf e}_s}=0$ for $s\gg 0$.

\begin{prop} \label{action-z=action-B} For each $t\geq 1$ and $\bfi\in{\scr B}_\infty$,
$$B^+_t(\wedge\oz_\bfi)=v^t \bdz_t^+\cdot(\wedge\oz_\bfi)\;\text{ and }\;
B^-_{t}(\wedge\oz_\bfi)=\bdz_t^-\cdot(\wedge\oz_\bfi).$$
\end{prop}

\begin{pf} As in \cite[(49)]{KMS}, for each $m\in\bbZ$, write
$$|m\ran=\oz_{m}\wedge\oz_{m-1}\wedge\oz_{m-2}\wedge\cdots \in \ooz^{\infty}/\Xi^\infty.$$
 Then $\bdz_t^+\cdot |m\ran=0$ and $K_\dz\cdot|m\ran=q|m\ran$. Write
$$\wedge\oz_\bfi=\oz_{i_1}\wedge\cdots\wedge \oz_{i_N}\wedge |N-m\ran.$$
 Applying \eqref{r-fold-z-pm} gives that
$$\aligned
{}&\bdz_t^+\cdot(\wedge\oz_\bfi)\\
=&\sum_{s=0}^{N} \underbrace{\oz_{i_1}\wedge\cdots\wedge \oz_{i_s}}_{s}\wedge \bdz^+_t\cdot\oz_{i_{s+1}}
\wedge \underbrace{K_{t\dz}\cdot\oz_{i_{s+2}}\wedge\cdots\wedge K_{t\dz}\cdot\oz_{i_N}}_{N-s-1}\wedge K_{t\dz}\cdot|N-m\ran\\
=&\sum_{s=0}^{N} v^{t}\underbrace{\oz_{i_1}\wedge\cdots\wedge \oz_{i_s}}_{s}\wedge \oz_{i_{s+1}+tn}
\wedge \underbrace{\oz_{i_{s+2}}\wedge\cdots\wedge \oz_{i_N}}_{N-s-1}\wedge \cdot|N-m\ran\\
=&v^{t} B^+_t(\wedge\oz_\bfi) \quad\text{(since $B_t^+(|N-m\ran)=0$),}
\endaligned$$
 that is, $B^+_t(\wedge\oz_\bfi)=v^t \bdz_t^+\cdot(\wedge\oz_\bfi)$. The second equality
 can be proved similarly.
\end{pf}

\begin{cor} Let $t,s\geq 1$. Then for each $\bfi\in{\scr B}_\infty$,
$${[\bdz_t^+,\bdz_s^-]}\cdot (\wedge\oz_\bfi)=\dz_{t,s}\, \frac{t(v^{2tn}-1)}{(v^t-v^{-t})^2}(K_{t\dz}-K_{-t\dz})\cdot(\wedge\oz_\bfi).$$
\end{cor}

\begin{pf} By \cite[Prop.~2.2~\&~2.6]{KMS} (with $q=v$),
$$[B_t^+,B_s^-]=\dz_{t,s}\frac{t(1-v^{2tn})}{1-v^{2n}}.$$
 This together with Proposition \ref{action-z=action-B} implies that for each $\bfi\in{\scr B}_\infty$,
$${[\bdz_t^+,\bdz_s^-]}\cdot (\wedge\oz_\bfi) =v^{t}[B_t^+,B_s^-]\dz_{t,s}\cdot (\wedge\oz_\bfi) =\dz_{t,s}\frac{tv^{t}(1-v^{2tn})}{1-v^{2n}}(\wedge\oz_\bfi).$$
 On the other hand,
$$\aligned
\dz_{t,s}\, \frac{t(v^{2tn}-1)}{(v^t-v^{-t})^2}(K_{t\dz}-K_{-t\dz})\cdot(\wedge\oz_\bfi)
&=\dz_{t,s}\, \frac{t(v^{2tn}-1)}{(v^t-v^{-t})^2}(v^t-v^{-t})(\wedge\oz_\bfi)\\
&=\dz_{t,s}\frac{tv^{t}(1-v^{2tn})}{1-v^{2n}}(\wedge\oz_\bfi).
\endaligned$$
 This gives the desired equality.
\end{pf}

In conclusion, $\bigwedge^\infty$ becomes a $\scrD(n)$-module which is obtained by
the restriction of the $\scrD(n)$-module structure on $\ooz^{\infty}/\Xi^\infty$ via the map $\kappa$.

\section{An isomorphism from $L(\llz_0)$ to $\bigwedge^\infty$}

In this section we show that the Fock space $\bigwedge^\infty$ as a $\scrD(n)$-module is isomorphic to
the basic representation $L(\llz_0)$ defined in Section 5. As an application, the decomposition of $L(\llz_0)$ in
Corollary \ref{decomp-hwm} induces the Kashiwara--Miwa--Stern decomposition of $\bigwedge^\infty$ in \cite{KMS}.

\begin{prop} \label{action-u-fkm}  For each $\fkm\in \fkM_n$, $\wt u^-_\fkm\cdot\bempty$ is a $\cZ$-linear combination
of those $\bmu$ satisfying $\fkm_\mu\leq_\dg \fkm$.
\end{prop}

\begin{pf} By \eqref{negative-action-Fock-space},
$$\wt u^-_\fkm\cdot\bempty=\sum_{\bfd}\big(\gz^-_\bfd(\wt u^-_\fkm)K_{-\bfd'}\big)\cdot\bempty,\;
\text{ where $\bfd'=\sum_{i}\big(\sum_{j<i,\,\bar j=\bar i}d_j\big)\vez_i$.}$$
 Since $K_i\cdot\bempty=v^{\dz_{i,0}}\bempty$ for $i\in\bbZ$, it follows
that $K_{-\bfd'}\cdot\bempty=v^{-\sum_{j<0,\,\bar j=\bar 0}d_j}\bempty$. By Proposition \ref{images-under-gamma},
$$\gz^-_\bfd(\wt u^-_\fkm)\in \sum_{\frak z}\cZ \wt u^-_{\frak z},$$
 where the sum is taken over ${\frak z}\in\fkM_\infty$ with $\scrF({\frak z})\leq^\infty_\dg \fkm$. Further,
by Proposition \ref{action-module-A-infty}(1),
$$\wt u^-_{\frak z}\cdot\bempty\in\cZ\bmu$$
 for some $\mu\in\Pi$ with $\fkm^\infty_\mu\leq^\infty_\dg \frak z$. This implies that
$$\fkm_\mu=\scrF(\fkm_\mu^\infty)\leq_\dg\scrF({\frak z})\leq_\dg\fkm.$$
 This finishes the proof.
\end{pf}

For each $\bfd=(d_i)\in\bbN I_\infty$, set
$$\sz(\bfd)=-\sum_{i<0,\,\bar i=\bar 0} d_i.$$
For $\lz\in\Pi$, we write $\sz(\lz)=\sz(\udim M(\fkm_\lz^\infty))$. The following
result was proved in \cite[9.2\,\&\,10.1]{VV}. We provide here a direct proof for completeness.

\begin{cor} \label{action-u-lz} For each $\lz\in\Pi$,
$$\wt u^-_{\fkm_\lz}\cdot\bempty\in \blz+\sum_{\mu\lhd\lz}\cZ \bmu.$$
In particular, the $\scrD(n)$-module $\bigwedge^\infty$ is generated by $\bempty$ and the set
$$\{b^-_{\fkm_\lz} \bempty\mid \lz\in\Pi\}$$
 is a basis of $\bigwedge^\infty$.
\end{cor}

\begin{pf} Applying Corollary \ref{image-m-lambda} gives that
$$\aligned
\wt u^-_{\fkm_\lz}\cdot\bempty&=\sum_{\bfd}\big(\gz^-_\bfd(\wt u^-_{\fkm_\lz})K_{-\bfd'}\big)\cdot\bempty
=\sum_{\bfd}v^{\sz(\bfd)} \gz^-_\bfd(\wt u^-_{\fkm_\lz})\cdot\bempty\\
&=\sum_{r\in\bbZ}v^{\thz(\lz)+\sz(\lz)}\wt u^-_{\tau^{rm}(\fkm^\infty_\lz)}\cdot\bempty +
\sum_{{\frak z}\in\fkM_\infty,\,\scrF({\frak z})<_\dg \fkm_\lz} f_{\lz,\frak z}\, \wt u_{\frak z}^-\cdot\bempty,
\endaligned$$
 where $f_{\lz,\frak z}\in\cZ$. By Proposition \ref{action-module-A-infty} and its proof,
$$\wt u^-_{\fkm^\infty_\lz}\cdot\bempty=\blz\;\text{ and }\;\wt u^-_{\tau^{rm}(\fkm^\infty_\lz)}\cdot\bempty=0\;
\text{for $r>0$}.$$
 Furthermore, for each $r<0$, $\wt u^-_{\tau^{rm}(\fkm^\infty_\lz)}\cdot\bempty\in\cZ\bmu$
 such that $\fkm_\mu^\infty\leq_\dg^\infty\tau^{rm}(\fkm^\infty_\lz)$.
Then $\fkm_\mu=\scrF(\fkm_\mu^\infty)\leq_\dg\scrF(\tau^{rm}(\fkm^\infty_\lz))=\fkm_\lz$, which implies
that $\mu\unlhd\lz$. Since $M(\tau^{rm}(\fkm^\infty_\lz))$ does not have a composition factor isomorphic to
$S_{\lz_1-1}$, $\mu$ does not contain a box with color $\lz_1-1$. Thus, $\mu\not=\lz$ and $\mu\lhd\lz$.

Finally, by Proposition \ref{action-u-fkm}, for each ${\frak z}\in\fkM_\infty$ with $\scrF({\frak z})<_\dg \fkm_\lz$,
$\wt u_{\frak z}^-\cdot\bempty$ is a $\cZ$-linear combination
of $\bmu$ satisfying $\fkm_\mu\leq_\dg \scrF(\frak z)$. Thus, $\fkm_\mu\leq_\dg\scrF({\frak z})<_\dg\fkm_\lz$,
which by Lemma \ref{order-coincidence} implies that $\mu\lhd\lz$. Hence, each $\wt u_{\frak z}^-\cdot\bempty$
is a $\cZ$-linear combination of $\bmu$ with $\mu\lhd\lz$. Consequently,
 $$\wt u^-_{\fkm_\lz}\cdot\bempty\in v^{\thz(\lz)+\sz(\lz)}\blz+\sum_{\mu\lhd\lz}\cZ \bmu.$$
Therefore, it remains to show that
$$\thz(\lz)+\sz(\lz)=0.$$
 Write $\lz=(\lz_1,\ldots,\lz_m)$ with $\lz_1\geq\cdots\geq \lz_m\geq 1$ and set $|\lz|=\sum_{s=1}^m\lz_s$.
We proceed induction on $|\lz|$ to show that $\thz(\lz)+\sz(\lz)=0$.
By the definition,
$$\thz(\lz)=\sum_{s<t}\kappa(\bfd_s,\bfd_t)-\sum_{s=1}^\ell h(\bfd_s),$$
 where $\ell=\lz_1$ is the Loewy length of $M=M(\fkm_\lz^\infty)$ and $S_{\bfd_s}\cong \rad^{s-1}M/\rad^s M$
for $1\leq s\leq\ell$. Let $1\leq t\leq m$ be such that
$\lz_1=\cdots =\lz_t>\lz_{t+1}$ and define
$$\lz'=(\lz_1,\cdots,\lz_{t-1},\lz_t-1,\lz_{t+1},\lz_m).$$
 Then $|\lz'|=|\lz|-1$. By the induction hypothesis, we have
$\thz(\lz')+\sz(\lz')=0$.

 For each $1\leq s\leq \ell$, let $\bfd_s'\in\bbN I_\infty$ be defined by setting
$S_{\bfd'_i}\cong \rad^{s-1}M'/\rad^s M'$, where $M'=M(\fkm_{\lz'}^\infty)$.
Then
$$\bfd_\ell'=\bfd_\ell-\vez_{\ell-t}\;\text{ and }\;\bfd'_s=\bfd_s\;\text{ for $1\leq s<\ell$.}$$
 This implies that
$$\sum_{s=1}^\ell h(\bfd_s)-\sum_{s=1}^\ell h(\bfd'_s)=h(\bfd_\ell)-h(\bfd_\ell')=-\dz_{\bar t,\bar 1}\;\text{ and }$$
$$\sum_{s<t}\kappa(\bfd_s,\bfd_t)-\sum_{s<t}\kappa(\bfd'_s,\bfd'_t)=\sum_{1\leq s<\ell}\kappa(\bfd_s,\vez_{\ell-t}).$$
 Hence,
$$\thz(\lz)-\thz(\lz')=\sum_{1\leq s<\ell}\kappa(\bfd_s,\vez_{\ell-t})+\dz_{\bar t,\bar 1}.$$
 On the other hand, $\sz(\lz)=\sz(\lz')-1$ if $\ell-t<0$ and $\bar\ell=\bar t$, and $\sz(\lz)=\sz(\lz')$ otherwise.
A direct calculation shows that if $\ell-t\geq 0$, then
$$\sum_{1\leq s<\ell}\kappa(\bfd_s,\vez_{\ell-t})=-\dz_{\bar t,\bar 1},$$
 and if $\ell-t< 0$, then
$$\sum_{1\leq s<\ell}\kappa(\bfd_s,\vez_{\ell-t})=\begin{cases} \dz_{\bar\ell,\bar t}-1, &\text{if $\bar t=\bar 1$};\\
                        \dz_{\bar\ell,\bar t}, &\text{if $\bar t\not=\bar 1$.}\end{cases}$$
We conclude that in all cases,
$$\thz(\lz)+\sz(\lz)=\thz(\lz')+\sz(\lz')=0.$$
\end{pf}

%
%

By the definition, for each $i\in I_n=\bbZ/n\bbZ$,
$$K_i\bempty=v^{\dz_{i,0}} \bempty.$$
 This together with the corollary above implies that $\bigwedge^\infty$ is a highest weight
$\scrD(n)$-module of highest weight $\llz_0$. Consequently, there is a unique surjective $\scrD(n)$-module homomorphism
$$\vphi: \scrD(n)^-=M(\llz_0)\lra \mbox{$\bigwedge^\infty$},\;\eta_{\llz_0}\lmto \bempty.$$

\begin{thm} \label{iso-brep-focks} The homomorphism $\vphi$ induces an isomorphism of $\scrD(n)$-modules
$$\bar\vphi:L(\llz_0)\lra \mbox{$\bigwedge^\infty$}.$$
\end{thm}

\begin{pf}  By definition, we have
$$F_i\cdot\bempty=0\;\text{ for $i\in I_n\backslash\{0\}$ and }\; F_0^2\cdot\bempty=0.$$
 Therefore, $\vphi$ induces a surjective homomorphism
$$\bar\vphi: L(\llz_0)=\scrD(n)^-/\big(\sum_{i\in I_n}\scrD(n)^- F_i^{\llz_0(h_i)+1}\big)\lra \mbox{$\bigwedge^\infty$}.$$
 Since $L(\llz_0)$ is simple, we conclude that $\bar\vphi$ is an isomorphism.
\end{pf}

Combining the theorem with Corollary \ref{decomp-hwm} gives the decomposition of $\bigwedge^\infty$
obtained by Kashiwara, Miwa and Stern in \cite[Prop.~2.3]{KMS}.

\begin{cor} As a $\bfU'_v({\widehat{\frak{sl}}}_n)$-module, $\bigwedge^\infty$ has a decomposition
$$\mbox{$\bigwedge^\infty$}|_{\bfU'_v(\widehat{\frak{sl}}_n)}\cong \bigoplus_{m\geq 0} L_0(\llz_0-m\dz^*)^{\oplus p(m)}.$$
\end{cor}

\section{The canonical basis for $\bigwedge^\infty$}

In this section we show that the canonical basis of $\bigwedge^\infty$
defined in \cite{LT} can be constructed by using the monomial basis of the Ringel--Hall
algebra of $\dt_n$ given in \cite{DDX}. We also interpret the ``ladder method'' in \cite{LLT}
in terms of generic extensions defined in Section 2.

Recall that there is a bar-involution $a\mapsto\iota(a)=\oline{a}$ on $\scrD(n)^-$ which
takes $\oline{v}\mapsto v^{-1}$ and fixes all $\tu_\az^-$ for $\az\in \bbN I_n$.
Then it induces a semilinear involution on the basic representation $L(\llz_0)$ by setting
$$\oline{a\eta_{\llz_0}}=\oline{a}\eta_{\llz_0}\;\text{ for all $a\in \scrD(n)^-$.}$$
 On the other hand, by \cite{LT}, there is a semilinear involution $x\mapsto \oline x$ on $\bigwedge^\infty$
which, by \cite{VV}, satisfies
\begin{itemize}
\item[(i)] $\oline{|\emptyset\ran}=|\emptyset\ran$,
\item[(ii)] $\oline {a x}=\oline{a}\,\oline{x}$ for all $a\in\scrD(n)^-$ and $x\in\bigwedge^\infty$.
\end{itemize}
Therefore, the isomorphism $L(\llz_0)\ra\bigwedge^\infty$ given in Theorem \ref{iso-brep-focks}
is compatible with the bar-involutions.

It is proved in \cite[Th.~3.3]{LT} that for each $\lz\in\Pi$,
\begin{equation}\label{bar-lz-presentation}
\oline{|\lz\ran}=|\lz\ran+\sum_{\mu\lhd\lz}a_{\mu,\lz}|\mu\ran,\;\text{ where $a_{\mu,\lz}\in\cZ$.}
\end{equation}
 Then applying the standard linear algebra method to the basis $\{|\lz\ran\mid\lz\in\Pi\}$
in \cite{L90} (or see \cite{Du} for more details) gives rise to an ``IC basis"
$\{b_\lz\mid \lz\in\Pi\}$ which is characterized by
$$\oline{b_\lz}=b_\lz\;\text{ and }\; b_\lz\in |\lz\ran+\sum_{\mu \lhd \lz}v^{-1}\bbZ[v^{-1}]|\mu\ran,$$
 The basis $\{b_\lz\mid\lz\in\Pi\}$ is called the {\it canonical basis} of $\bigwedge^\infty$.
In other words, the basis elements $b_\lz$ are uniquely determined by the polynomials
$a_{\mu,\lz}$.

\begin{rem} Varagnolo and Vasserot \cite{VV} have conjectured that
$$b^-_{\fkm_\lz}\cdot|\emptyset\ran=b_\lz\;\text{ for each $\lz\in\Pi$.}$$
 This conjecture was proved by Schiffmann \cite{Sch00}.
\end{rem}

In the following we provide a way to deduce \eqref{bar-lz-presentation} by using
the monomial basis of the Ringel--Hall algebra of $\dt_n$ given in \cite{DDX}.
As in \cite[Sect.~3]{DDX}, set
$$I^e= I_n\cup\{\text{all sincere vectors in $\bbN I_n$}\}$$
and consider the set $\Sigma$ of all words on the alphabet $I^e$. Since $\scrD(n)^-$ is isomorphic to the opposite
Ringel--Hall algebra of $\dt_n$, we define
$$M\ast' N=N\ast M.$$
 This gives the map
$$\wp^{\rm op}:\Sigma\lra \fkM, \;w=\bfa_1\bfa_2\cdots \bfa_t\longmapsto S_{\bfa_1}\ast' S_{\bfa_2}\ast'\cdots\ast' S_{\bfa_t}.$$
  By \cite[Sect.~9]{DDX}, for each $\fkm\in\fkM$, there is a distinguished
word $w_\fkm\in({\wp^{\rm op}})^{-1}(\fkm)$ which defines a monomial $m^{(w_\fkm)}$ on $\wt u_\bfa^-$
with $\bfa\in \wt I$ such that
$$m^{(w_\fkm)}=\tu_\fkm^-+\sum_{\fkp<_\dg\fkm}\thz_{\fkp,\fkm}\tu_\fkp^-\;\text{ for some $\thz_{\fkm,\fkp}\in \cZ$;}$$
 see \cite[(9.1.1)]{DDX}. If $\fkm=\fkm_\lz$ for some $\lz\in\Pi$, we simply write $w_{\fkm_\lz}=w_\lz$.
Thus,
\begin{equation}\label{mon-PBW-0}
m^{(w_\lz)}=\tu_{\fkm_\lz}^-+\sum_{\fkp<_\dg\fkm_\lz}\thz_{\fkp, \fkm_\lz}\tu_\fkp^-.
\end{equation}
 This together with Proposition \ref{action-u-fkm} and Corollary \ref{action-u-lz} implies that
\begin{equation}\label{mon-PBW}
m^{(w_\lz)}|\emptyset\ran=|\lz\ran+\sum_{\mu\lhd\lz}\tau_{\mu,\lz}|\mu\ran,
\end{equation}
 where $\tau_{\mu,\lz}\in\cZ$. Since the monomials $m^{(w_\lz)}$ are bar-invariant, we deduce that
for each $\lz\in\Pi$,
$$\oline{|\lz\ran}=|\lz\ran+\sum_{\mu\lhd\lz}a'_{\mu,\lz}|\mu\ran\;\text{ for some $a'_{\mu,\lz}\in\cZ$.}$$
 Comparing with \eqref{bar-lz-presentation} gives that
$$a_{\mu,\lz}=a'_{\mu,\lz}\;\text{ for all $\mu\lhd\lz$.}$$

In case $\lz$ is $n$-regular, then $\fkm_\lz$ is aperiodic and the word $w_\lz$ can be chosen
in $\ooz$, the subset of all words on the alphabet $I_n=\bbZ/n\bbZ$; see \cite[Sect.~4]{DDX}. In other words, $m^{(w_\lz)}$
is a monomial of the divided powers $(u_i^-)^{(t)}=F_i^{(t)}$ for $i\in I_n$ and $t\geq 1$. We now
interpret the ``ladder method'' in \cite[Sect.~6]{LLT} in terms of the generic extension map. Let
$\lz=(\lz_1,\ldots,\lz_t)\in\Pi$ be $n$-regular. Recall the corresponding nilpotent representation
$$M(\fkm_\lz)=\bigoplus_{a=1}^t S_{1-a}[\lz_a],$$
 where $1-a$ is viewed as an element in $I_n$. Take $1\leq s\leq t$ with
$\lz_1=\cdots=\lz_s>\lz_{s+1}$ ($\lz_{t+1}=0$ by convention) and let $k\geq 0$ be maximal such that
$$\lz_{s+l(n-1)+1}=\cdots=\lz_{s+(l+1)(n-1)}\text{ and }\;
\lz_{s+l(n-1)}=\lz_{s+l(n-1)+1}+1\;\text{for $0\leq l\leq k-1$}.$$
 Let $i_1\in I$ be such that $\soc(S_{1-s}[\lz_s])=S_{i_1}$. Then for each $a=s+l(n-1)$
with $0\leq l\leq k$,
$$\soc(S_{1-a}[\lz_a])=S_{i_1}.$$
Define $\mu=(\mu_1,\ldots,\mu_t)\in\Pi$ by setting
$$\mu_a=\begin{cases} \lz_a-1, &\text{if $a=s+l(n-1)$ for some $0\leq l\leq k$};\\
                        \lz_a, &\text{otherwise.}\end{cases}$$
 It is easy to see from the construction that $\mu$ is again $n$-regular. Moreover, by applying
an argument similar to that in the proof of \cite[Prop.~3.7]{DD1},
$$(k+1)S_{i_1}\ast' M(\fkm_\mu)=M(\fkm_\mu)\ast (k+1)S_{i_1}= M(\fkm_\lz).$$
 Repeating the above process, we finally obtain a sequence $i_1,\ldots,i_d$ in $I_n$ and positive
 integers $k_1=k+1,\ldots,k_d$ such that
$$(k_1S_{i_1})\ast'\cdots\ast' (k_dS_{i_d})=M(\fkm_\lz).$$
 In other word, the word $w_\lz:=i_1^{k_1}\cdots i_d^{k_d}$ lies in $(\wp^{\rm op})^{-1}(\fkm_\lz)$.
It can be also checked that the word $w_\lz$ is distinguished. Thus, the corresponding monomial
$$m^{(w_\lz)}=(u_{i_1}^-)^{(k_1)}\cdots (u_{i_d}^-)^{(k_d)}=F_{i_1}^{(k_1)}\cdots F_{i_d}^{(k_d)}$$
gives rise to the equality \eqref{mon-PBW} for the element $m^{(w_\lz)}|\emptyset\ran$.
We remark that $m^{(w_\lz)}|\emptyset\ran$ coincides with the element $A(\lz)$ constructed
in \cite[(8)]{LLT} by using the ``ladder method'' of James and Kerber \cite{JK}.

\end{document}